\documentclass{article}
\usepackage[english]{babel}
\usepackage{bookmark}
\usepackage[letterpaper,top=2cm,bottom=2cm,left=3cm,right=3cm,marginparwidth=1.75cm]{geometry}

\usepackage{amsmath}
\usepackage[table,xcdraw]{xcolor}
\usepackage{amsthm}
\usepackage{amssymb}
\usepackage{graphicx}
\usepackage{amsfonts}
\usepackage{tikz}
\usepackage{float}
\usepackage{longtable}
\usepackage{mathtools}
\usepackage{enumerate}
\usepackage{dsfont}
\usepackage{bbm}
\usepackage{bm}
\usepackage{booktabs}
\usepackage{nicematrix}
\usepackage{multirow}
\usepackage{prettyref}
\usepackage{mathtools}
\usepackage{listings}
\usepackage{xcolor}
\lstdefinestyle{RStyle} {
    language=R,
    style=generalStyle  
}
\definecolor{mygreen}{rgb}{0,0.6,0}
\definecolor{mygray}{rgb}{0.5,0.5,0.5}
\definecolor{mymauve}{rgb}{0.58,0,0.82}
\usepackage{algorithm}
\usepackage{algorithmic}
\usepackage{soul}
\usepackage{subfigure}
\usepackage{appendix}

\usetikzlibrary{arrows,intersections}

\newtheorem{definition}{Definition}[section]

\title{\vspace{-0.9cm} Achieving Fairness and Accuracy in Regressive Property Taxation}
\author{ Ozan Candogan \thanks{The University of Chicago, Booth School of Business (ozan.candogan@chicagobooth.edu).}\and Feiyu Han \thanks{The University of Chicago, Booth School of Business  (feiyu.han@chicagobooth.edu).}\and Haihao Lu\thanks{The University of Chicago, Booth School of Business (haihao.lu@chicagobooth.edu).} }

\begin{document}
\date{}
\maketitle

\begin{abstract}
% \textit{Problem Definition:}
Regressivity in property taxation, or the disproportionate overassessment of lower-valued properties compared to higher-valued ones, results in an unfair taxation burden for Americans living in poverty. To address regressivity and enhance both the accuracy and fairness of property assessments, we introduce a scalable property valuation model called the $K$-segment model. Our study formulates a mathematical framework for the $K$-segment model, which divides a single model into $K$ segments and employs submodels for each segment. Smoothing methods are incorporated to balance and smooth the multiple submodels within the overall model. To assess the fairness of our proposed model, we introduce two innovative fairness measures for property evaluation and taxation, focusing on group-level fairness and extreme sales price portions where unfairness typically arises. Compared to the model employed currently in practice, our study demonstrates that the $K$-segment model effectively improves fairness based on the proposed measures. Furthermore, we investigate the accuracy--fairness trade-off in property assessments and illustrate how the $K$-segment model balances high accuracy with fairness for all properties. Our work uncovers the practical impacts of the $K$-segment models in addressing regressivity in property taxation, offering a tangible solution for policymakers and property owners. By implementing this model, we pave the way for a fairer taxation system, ensuring a more equitable distribution of tax burdens. 
\end{abstract}

\section{Introduction} \label{back}

Property assessments are a critical aspect of local government operations, as they determine the value of a property for tax purposes. These assessments are often used in conjunction with local tax rates to calculate the amount of property tax each property owner is required to pay, see, e.g., \cite{cookcountyassessor}.
Therefore, the accuracy of these assessments is vital to ensure fair taxation.

In some cases, such as in Chicago \cite{berry2021reassessing,mcmillen2020assessment}, inaccurate property assessments have been argued to result in what is termed ``regressive taxation.'' This occurs when properties at the lower end of the market are assessed at values higher than their actual worth (often measured by sales price). This leads to a higher effective tax rate for these properties, disproportionately affecting owners of less valuable properties.
The impact of this regressive taxation can be substantial. For instance, owners of lower-valued properties might face effective tax rates up to 50\% higher than those paid by owners of higher-valued properties. This disparity exists despite both groups having access to the same amenities funded by property taxes and being subject to the same statutory property tax rates \cite{amornsiripanitch2020residential}.
This issue of regressive taxation due to inaccurate property assessments is not unique to Chicago but is a pattern that is argued to be prevalent throughout the United States. The limitations of the data and assessment methods used are often cited as reasons for these inaccuracies; see, e.g., \cite{berry2021reassessing}.

% Conventionally, assessors determine the value of properties using the sale prices of similar properties. However, properties with similar observable features can vary significantly in value. %Assessors may not be sensitive to intangible factors such as local market conditions and the neighborhood's unique characteristics. 
% Consequently, owners of lower-valued properties can face up to 50\% higher effective tax rates than owners of higher-valued properties, despite having access to the same property tax-funded amenities and statutory property tax rates \cite{amornsiripanitch2020residential}. This disproportionately penalizes low-income individuals, underscoring the urgent need for policymakers and property owners to address property assessment regressivity.

This paper
takes the position that regressive taxation and assessment errors are artifacts of the machine learning algorithms used to assess property values. It investigates how these
methods can be modified to alleviate regressive taxation. 
Considering the complexity of machine learning pipelines used in practice and the difficulty of influencing them, the paper focuses on pinpointing minimal adjustments to current pipelines to achieve this objective.

More precisely, our analysis is based on the 2022 residential reassessment model, hereafter the \textbf{original model}, of the Cook County Assessor's Office (CCAO). This model uses the attributes of the properties and the sales price data collected by the CCAO to assess the values of all the properties in Cook County. %This model will be referred to as the `original model' throughout this paper. 
We detail this model and provide evidence for its regressivity in Section \ref{eva}.

Regressive taxation is often viewed as unfair.
Although fair classification has been a significant focus of recent research, the concept of fairness in regression models has been less extensively explored in applied work related to property taxation, which requires novel context-specific notions of fairness.
 Our first contribution is to introduce two novel fairness notions
 to quantify the amount of unfairness/regressivity
created by different assessment models.
The first notion adapts existing group fairness concepts from the literature to our context by examining the \emph{ratios between assessed values and actual property sales values}. 
The second notion, which we call deviation-weighted fairness, emphasizes penalizing deviations from a ratio of one. Importantly, it uses weighting functions that impose heavier penalties on properties with exceptionally low or high values. 
Intuitively, such weighting functions are particularly relevant in addressing regressive taxation issues, as these issues tend to be more pronounced in properties with extreme values.

The original model causing regressivity is based on the Light Gradient Boosting Machine (LightGBM). Our second contribution involves proposing a simple modification to this model, leading to the development of our $K$-segment model.
This model divides the data points used by the original model into
$K$ segments, applying the original model independently to each segment. These segments are determined by using the prior year's assessments as input and partitioning the property set accordingly. For example, in a 5-segment model ($K=5$), the lowest 20\% of the properties are grouped into the first segment, the next 20\% into the second, and so on. Intuitively, this approach ensures that the model is separately trained on properties with inherently similar values, thereby potentially reducing the extent of regressive taxation within each segment.

However, naively applying this segmentation idea is problematic. 
This is because, the resulting assessment model is inherently ``discontinuous.'' More precisely, properties that are to the left/right of a cutoff assessment value that defines the segments are assessed using different models. Hence, in the vicinity of these cutoffs, small changes to assessed property values in preceeding years can result in large assessment differences in subsequent years.

Our third contribution is to introduce a set of smoothing methods that eliminate such discontinuous outcomes. These methods
use an exponential weighting approach (or variants thereof) to
``blend'' segments, so that some properties that belong to one segment in the approach described above are transferred to other segments (where the likelihood of such transfers decays exponentially for non-adjacent segments).

Finally, we use our model to produce assessments for Cook County, and demonstrate that
compared to the original model, our $K$-segment model yields
assessments that ensure that 
lower- (higher-)valued properties (according to their sales prices) are assessed lower (higher). Thus our approach mitigates regressivity.
We also quantify the improvement using the fairness notions described earlier.
Perhaps more surprisingly, we show that the resulting assessments are also more accurate, and that the $R^2$ of our $K$-segment model improves over the original model. Thus, in a sense, our method provides a ``Pareto-improvement'' on the original model.

All in all, our work provides a concrete recipe for alleviating regressive property taxation by incorporating a minimal change to the existing data science pipeline, thus holding the promise of improving property assessments in practice.
More broadly, our work represents a novel contribution to the field of fair regression in data science and machine learning and has potential implications for fair evaluation practices beyond the housing market.

\textbf{Outline.}
The paper is organized as follows. 
Section~\ref{fair} provides two novel measures of fairness for measuring regressivity. 
Section~\ref{model} formally introduces the $K$-segment model, along with its smoothing methods. 
Section~\ref{eva} evaluates the $K$-segment model in terms of smoothness, accuracy, and fairness, while exploring the trade-off between accuracy and the two fairness measures presented in Section~\ref{fair}. Section~\ref{con} concludes and suggests future directions. 
Additional results on the unsmoothed model and other parameter choices for the smoothing method can be found in the appendix.

\textbf{Related Work.}
 Several theories have been proposed that explain the regressivity of property value assessments \cite{paglin1972equity,weber2010ask,amornsiripanitch2020residential}. Much remains unexplained, however; data and machine learning methods could be instrumental in addressing this issue. The measure of regressivity in assessments primarily depends on the indicator used. The sales-to-assessment (StA) ratio is the most common indicator. Nonetheless, other measures can also gauge regressivity. For instance, \cite{cornia2005property} employed the coefficient of dispersion, representing the average percentage deviation of assessment ratios from the median ratio. \cite{mcmillen2020assessment} utilized the price-related differential, calculated as the median assessment ratio for lower-valued properties divided by the same ratio for higher-valued properties.

Fairness in machine learning is an active research area.\footnote{The discussions have predominantly focused on fair classification.}
A seminal work on fair regression \cite{calders2013controlling} employed mean constraints in linear regression models to ensure equitable treatment across different groups, thereby mitigating bias and enhancing fairness.
\cite{berk2017convex} proposed a set of versatile fairness regularizers for linear and logistic regression problems, while \cite{komiyama2018nonconvex} considered imposing a fairness constraint, not a regularizer.
\cite{agarwal2019fair} suggested capturing fairness by minimizing the expected loss of predictions with real value under conditions of statistical parity or bounded group loss. Meanwhile, \cite{diana2021minimax} proposed a fairness measurement by worst-case outcomes across groups, using a minimax optimization method.
A myriad of other existing methods include decision trees to ensure disparate impact and treatment in regression tasks \cite{raff2018fair,aghaei2019learning}, a method for considering multiple sensitive attributes simultaneously \cite{perez2017fair}, and algorithms learning fair predictors \cite{agarwal2019fair}. Other notable methodologies include an automated procedure addressing fairness in decision-making when sensitive features are involved \cite{pmlr-v97-nabi19a}, and a group fairness enforcement mechanism as a constraint on kernel regression \cite{fitzsimons2019general}. \cite{shah2022selective} focused on performance disparities between subgroups in selective regression, proposing neural network-based algorithms to mitigate disparities.
Several studies have looked into ``pairing individuals" in the context of fairness in regression and ranking models, involving cross-pairs, i.e., pairs of instances from different groups \cite{berk2017convex, lahoti2019operationalizing, bird2019fairness, kleindessner2022pairwise}. This approach is relevant to our group fairness concept, as outlined in Definition~\ref{grp}.

Striving for fairness while maintaining accuracy is a two-pronged optimization challenge in decision-making.
\cite{menon2018cost} proposed a method in which fairness is considered as a constraint on subgroup performance within the model, employing Lagrangian relaxation. \cite{zafar2017fairness} developed a flexible convex optimization framework that tries to balance accuracy and fairness by minimizing the accuracy loss function, while subjecting it to fairness constraints. 
\cite{zafar2019fairness} introduced a decision boundary covariance constraint, furthering the effort to harmonize accuracy and fairness.
\cite{dutta2020there} approached fair classification from the perspective of mismatched hypothesis testing, and \cite{pmlr-v80-agarwal18a} introduced reductions to calibrate the trade-off between accuracy and fairness.
Our work uses the concept of Pareto frontier to identify and select the best fairness solutions, i.e., those that are not dominated by other fairness solutions, without compromising accuracy (see Section~\ref{se3}). Further methods for identifying and selecting optimal solutions on the Pareto frontier have been explored in \cite{lotov2004interactive,valdivia2021fair}.

\section{Preliminaries: Quantification of Fairness}\label{fair}
% Two different aspects: group fairness and end sensitive fairness.
% The first aspect, group fairness, is concerned with ensuring that the model's predictions are similar for different groups on average. The model achieves this by dividing the data into groups and ensuring that each group is treated fairly.
% The second aspect, end sensitive fairness, is concerned with ensuring that the model is fair across the entire range of data, including the extremes. The model achieves this by using weighting functions that give more weight to extreme values, which allows the model to capture how unfair it is compared to an ideal situation.

We start by developing a mathematical measure of the extent of regressivity. One common method of diagnosis involves examining the sales-to-assessment (StA) ratios of properties with varying sales prices. The StA ratio of a property is determined by dividing its assessed value by its sales price. Thus, the ideal valuation of the StA ratio is 1.

The StA ratio may not be enough to determine fairness in terms of regressivity.
First, the StA ratio is confined to providing a singular value for each property. This design inherently neglects the need for a relational or comparative analysis between properties, thereby restricting the ability to evaluate how assessments are trending in relation to different price levels.
Furthermore, different properties sold at the same price often have varying StA ratios due to factors such as property characteristics and market conditions. This inconsistency hinders the ability to identify a common indicator for a specific sales price. Averaging the StA ratios for properties with the same sales price may seem a solution, but the presence of extreme StA ratios, either significantly higher or lower than the rest, can distort this average, leading to a misleading representation of that particular price range.
Since assessing regressivity involves understanding the pattern of overestimation or underestimation of properties across different price ranges, the StA ratio's limitations in comparing property' assessments within the measure, and its failure to represent specific prices, make it an inadequate measure. Therefore, more comprehensive measures are required to analyze and address regressivity.

To establish a more reliable framework for property valuation, we employ a multifaceted assessment model, in which various factors collectively determine the property's valuation.
% Besides, we particularly focus on $x_p\in\mathbb{R}_+$, which is the assessment value from the previous year's model. This value serves as an indicator of the property's perceived expense.
% We represent the quantile of $x_p$ as \( y \in\mathcal{Y}=[0,1] \). The introduction of the quantile \( y \) captures the relative price relationship between properties. It is computed as
% \[
% y = \frac{N(x_p)}{m_T},
% \]
% where \( N(x_p) \) indicates the count of properties whose price assessed in previous year's model less than or equal to \( x_p \), and $m_T$ represents the total number of properties in the training dataset. 
In the model's evaluation, the dataset consists of \( m \) properties with the known sales price $x\in\mathbb{R}_+$. We define the quantile for the sales price $x$ as $\Tilde{y}\in\mathcal{Y}=[0,1]:=N(x)/m$, where \( N(x) \) indicates the number of properties whose sales price is less than or equal to \( x \).
For each property, we obtain the predicted property valuation \( v\in \mathbb{R}_+ \) from the model, which allows us to calculate the StA ratio \( r \in \mathbb{R}_+=v/x\). This information forms the dataset
\[
D = \{(x_i, \Tilde{y}_i, v_i, r_i)\}_{i=1}^m,
\]
extracted from the underlying distribution.

Before we introduce our first fairness notion, we establish the concept of groups. We partition \(\mathcal{Y}\) into \(n\) groups \(G_1, \ldots, G_n\), where each group \(G_l\) contains \(m_l\) samples such that \(m = m_1 + \ldots + m_n\). 
Each group contains a distinct range of $\Tilde{y}$ values.
%, corresponding to the quantiles of assessed value of the model from the previous year.
Group $1$ encompasses the lowest range, and each subsequent group represents a progressively higher segment of sales prices, culminating in Group $n$, which contains the highest range. 
% Next, we define group fairness. 

\begin{definition}\label{grp}
    [Group Fairness]
   Let Group \(\sigma(i)\) designate the group to which \(\Tilde{y}_i\) belongs, and let \(m_{\sigma(i)}\) represent the number of samples in Group \(\sigma(i)\). 
   The group fairness measure \(F_{\text{grp}}\) is established by evaluating the differences between the StA ratios \(r_i\) and \(r_j\), associated with $\Tilde{y}_i$ and $\Tilde{y}_j$ from distinct groups, where \(\sigma(i)<\sigma(j)\):
    \[F_{\text{grp}}= -\sum_{i,j:\;\sigma(i)<\sigma(j)}\frac{1}{m_{\sigma(i)}m_{\sigma(j)}}(r_i-r_j)^+.\]
    % where $\sigma(i)$ designates the group to which data point $y_i$ belongs. 
\end{definition}

The group fairness measure, denoted by $F_{grp}$, quantifies the differences in the StA ratios between pairs of samples from distinct groups. Specifically, it calculates the positive differences\footnote{To calculate $F_{grp}$, we employ the positive part function, denoted by $(\cdot)^+$.} between the StA ratios of sample pairs $(\Tilde{y}_i, \Tilde{y}_j)$, where $\Tilde{y}_i$ belongs to a group with a lower index than $\Tilde{y}_j$'s group.
This method accounts for the phenomenon of regressivity, whereby overestimated values are more common in groups with lower indices, by primarily considering instances where the StA ratio \( r \) of a group with a smaller index exceeds ratios of groups with larger indices. When sample \( \Tilde{y}_i \) is part of group \( G_{\sigma(i)} \) with a lower index \( \sigma(i) \), and \( \Tilde{y}_j \) belongs to group \( G_{\sigma(j)} \) where \( \sigma(j) > \sigma(i) \), the fairness measure accounts for any surplus of \( r_i \) over \( r_j \). In the case of \( r_j > r_i \) arises, this surplus is excluded from the fairness measure, registering as \( 0 \).\footnote{In our work, we redefine group fairness in the context of regressivity issues, by employing the concept of cross-pairing from earlier studies such as \cite{berk2017convex, lahoti2019operationalizing}. This context-specific definition captures scenarios where groups with lower explanatory values have disproportionately higher target values compared to other groups with higher explanatory values.}

An alternative approach to defining the measure is to capture both overestimation and underestimation. This method incorporates weighting functions $\{(\Tilde{w}_1(\Tilde{y}_i),\Tilde{w}_2(\Tilde{y}_i))\}_{i=1}^m$ to penalize extreme quantiles more severely:
\begin{definition}\label{dyn}
[Deviation-Weighted Fairness]
We introduce weighting functions $\{(\Tilde{w}_1(\Tilde{y}_i),\Tilde{w}_2(\Tilde{y}_i))\}_{i=1}^m$ that are not reliant on the structure of the regression model, where $\Tilde{w}_1$ decreases and $\Tilde{w}_2$ increases with $\Tilde{y}\in\mathcal{Y}$. We define deviation-weighted fairness as follows:
% \[F_{dev}=-\Big[\int_0^1(r(y)-1)^+\cdot w_1(y) dy+\int_0^1(1-r(y))^+\cdot w_2(y)dy\Big],\]
\[F_{dev}=-\Bigg[\sum_{i=1}^m(r_i-1)^+\cdot \Tilde{w}_1(\Tilde{y}_i) + \sum_{i=1}^m(1-r_i)^+\cdot \Tilde{w}_2(\Tilde{y}_i)\Bigg].\]
 In the remainder of the paper, we set the weighting functions to ba $\Tilde{w}_1(\Tilde{y})=e^{-\alpha \Tilde{y}}, \Tilde{w}_2(\Tilde{y})=e^{-\alpha(1-\Tilde{y})}$, where $\alpha\geq 0$.
\end{definition}
 Deviation-weighted fairness, denoted by $F_{dev}$, is comprised of two weighting functions $\Tilde{w}_1$ and $\Tilde{w}_2$ that are utilized to assign greater importance to the extremes of the quantiles. These weighting functions accentuate overestimation in the lower ends and underestimation in the higher ends, respectively. By emphasizing the ends of the sales price range, where regressivity tends to be prominently observed in real-world situations, $F_{dev}$ can more effectively capture the degree of unfairness in a model relative to an ideal scenario, where a model achieves perfect assessment without any regressivity.

Overall, group fairness ($F_{grp}$) mainly assesses intergroup relations, focusing on differences in StA ratios between different groups, and penalizing larger disparities. It tries to ensure that the lower end is not
systematically disadvantaged relative to the higher end.
Deviation-weighted fairness ($F_{dev}$) emphasizes overall model conformity to an ideal baseline StA ratio of 1. By employing weighting functions that emphasize the extremes of the prediction distribution, $F_{dev}$ accounts for the frequently observed overestimation and underestimation at the respective ends in property assessment models. Together, these two fairness measures highlight different aspects of fairness.

In Section~\ref{se3}, we illustrate how the fairness measures $F_{grp}$ and $F_{dev}$ are used to assess fairness performance in various models using real data. 
The analysis demonstrates that they are both reliable indicators alongside accuracy measures in evaluating a model; due to their distinct fairness assessments, however, their values might not always coincide.

\section{\texorpdfstring{$K$}{K}-segment Model}\label{model}
In this section, we propose a novel method for obtaining ``fairer'' assessments.
Specifically, we introduce the $K$-segment model, which splits the data into $K$ distinct segments, and fits the baseline model to each segment separately. 
By appropriately tuning $K$, our approach achieves different trade-offs between model accuracy and fairness.
Naively splitting the data into $K$ quantiles and training a separate model on each segment results in ``discontinuity'' of the estimates in the vicinity of splitting points. To address this issue, we also develop different smoothing methods.
In the remainder of this section, we discuss both the unsmoothed model and the following versions of the $K$-segment model for which we use different smoothing methods:
\begin{itemize}
    \item $M_{\text{unsm},K}$ for the unsmoothed $K$-segment model;
    \item $M_{\text{q},K}$ for the quantile-based smoothing $K$-segment model;
    \item $M_{\text{m-s},K}$ for the midpoint score-based smoothing $K$-segment model;
    \item $M_{\text{d-s},K}$ for the distance score-based smoothing $K$-segment model.
\end{itemize}
We next formally define these models.

\subsection{Model Formulation}\label{modelfor}

% \subsubsection{From the Original to the \texorpdfstring{$K$}{K}-Segment Model}
In 2022, the CCAO implemented a comprehensive seven-stage pipeline for the assessment of residential properties in Cook County, ranging from raw data ingestion to the final output.
Within this pipeline, our primary focus is on the training and assessment stages.
The training stage employs the Light Gradient Boosting Machine (LightGBM) algorithm,\footnote{LightGBM is a gradient-boosting decision tree framework \cite{ke2017lightgbm}, frequently used in housing-specific machine learning models due to its high precision, speed, and comprehensive documentation.} iteratively updating model parameters to reduce errors and recognize patterns in sales data. 
The training is based on the known sales prices of similar and nearby properties, and takes into account factors such as property characteristics, location, environmental variables, and market trends. 
The assessment stage systematically predicts the market value of each unsold property across the county.\footnote{For readers seeking more in-depth information on the whole pipeline, additional details can be found in \cite{CCAO2022}.}
Hereafter, we refer to the model that combines these two stages as the \textbf{original model}.

In practice, it is observed that the original model leads to regressivity in assessments. To combat this and other modeling drawbacks, some post-modeling adjustments have been implemented (see \cite{CCAO2022}). However, despite these adjustments, the regressivity issue persists for properties sold in 2021 (as we demonstrate in Section~\ref{plot}).
%\footnote{Regressivity in property taxation can indeed be influenced in part by the nature of the data's distribution and the modeling method employed in assessments. This issue becomes pronounced toward the extremely low and comparatively high end of the price spectrum, see Figure~\ref{skew} in Appendix~\ref{app4}. }
% However, the complexity of LightGBM can pose interpretability challenges, and even with post-modeling adjustments to the predicted value\footnote{The original model employed post-modeling to rectify minor prediction errors, with the intention of ensuring uniform valuation for nearly identical properties and to reduce regressivity in assessment.}, it may still exhibit the same regressivity that most assessors and machine learning practitioners face, as we will show in Section~\ref{3.1}.
To address this challenge, we propose the $K$-segment model.
Specifically, in this model, we divide the quantile of assessed sales prices derived from last year's model,\footnote{The assessment value from the previous year's model, $x_p\in\mathbb{R}_+$, serves as an indicator for the property's perceived expense.
 We represent the quantile of $x_p$ as \( y \in\mathcal{Y}=[0,1] \). The introduction of the quantile \( y \) captures the relative price relationship between properties. It is computed as $y = N(x_p)/m_T$, where \( N(x_p) \) is the number of properties whose price assessed in the previous year's model is less than or equal to \( x_p \), and $m_T$ is the total number of properties in the training dataset.} $\mathcal{Y}$, into $K$ distinct segments and
apply the original model to each segment separately. Intuitively, doing so allows the model to capture the specific characteristics and patterns
of the data within each segment better, and
reduces the systematic errors impacting extremely low- and high-end properties.

We denote the segmentation parameters by \(\eta=(\eta_0,\eta_1,\ldots,\eta_K)\), where the thresholds divide the data into \(K\) distinct segments. Specifically, the \(k\)-th segment is defined as the interval \([\eta_{k-1}, \eta_k]\), with \(\eta_0=0\) and \(\eta_K=1\), and \(\eta_{k-1}<\eta_k\) for each segment.
Each segment corresponds to a unique submodel, \(S_k\), that is responsible for assessing value when the quantile $y$ falls within its respective \(k\)-th segment\footnote{In our implementation, properties whose
quantile equals some $\eta_k$ are assigned to the segment $k+1$; except in the case where $k=K$. This convention is used to break ties, and has negligible impact on our results. To ease notation, with some abuse, we used closed brackets to define each segment.} \([\eta_{k-1}, \eta_k]\). These submodels are trained on their segment-specific data, adhering to the training procedure of the original model.
Collectively, these submodels form the complete \(K\)-segment model, denoted by \(M(y)\).

The most basic form of this approach is the \textbf{Unsmoothed 
$K$-segment Model}, constructed using $K$ distinct and discontinuous submodels, each denoted by $S_k$, where the index $k$ ranges from 1 to $K$,
\[
M_{\text{unsm}}(y) =
S_k\big(y\big),  \text{if } y \in [\eta_{k-1},\eta_{k}],
\]
and a property with quantile  $y\in\mathcal{Y}$ is categorized into one of $K$ segments.
% \[
% M_{\text{unsm}}(y) =
% \begin{cases}
% S_1\big(y\big), & \text{if } y \in [\eta_{0},\eta_{1}], \\
% S_2\big(y\big), & \text{if } y \in [\eta_{1},\eta_{2}], \\
% \vdots \\
% S_K\big(y\big), & \text{if } y \in [\eta_{K-1},\eta_{K}].
% \end{cases}
% \]

% \subsubsection{Smoothing Techniques}

However, the unsmoothed model can create abrupt transitions at segment boundaries $\eta_k$ (see Appendix~\ref{app1}), leading to interpretability challenges. Properties near a segmentation boundary can be assigned different valuations despite having similar price percentiles, as different submodels are applied. Such inconsistency in valuation underscores the need for smoothing methods.

\begin{figure}[H]
    \centering
    \includegraphics[width=0.6\textwidth]{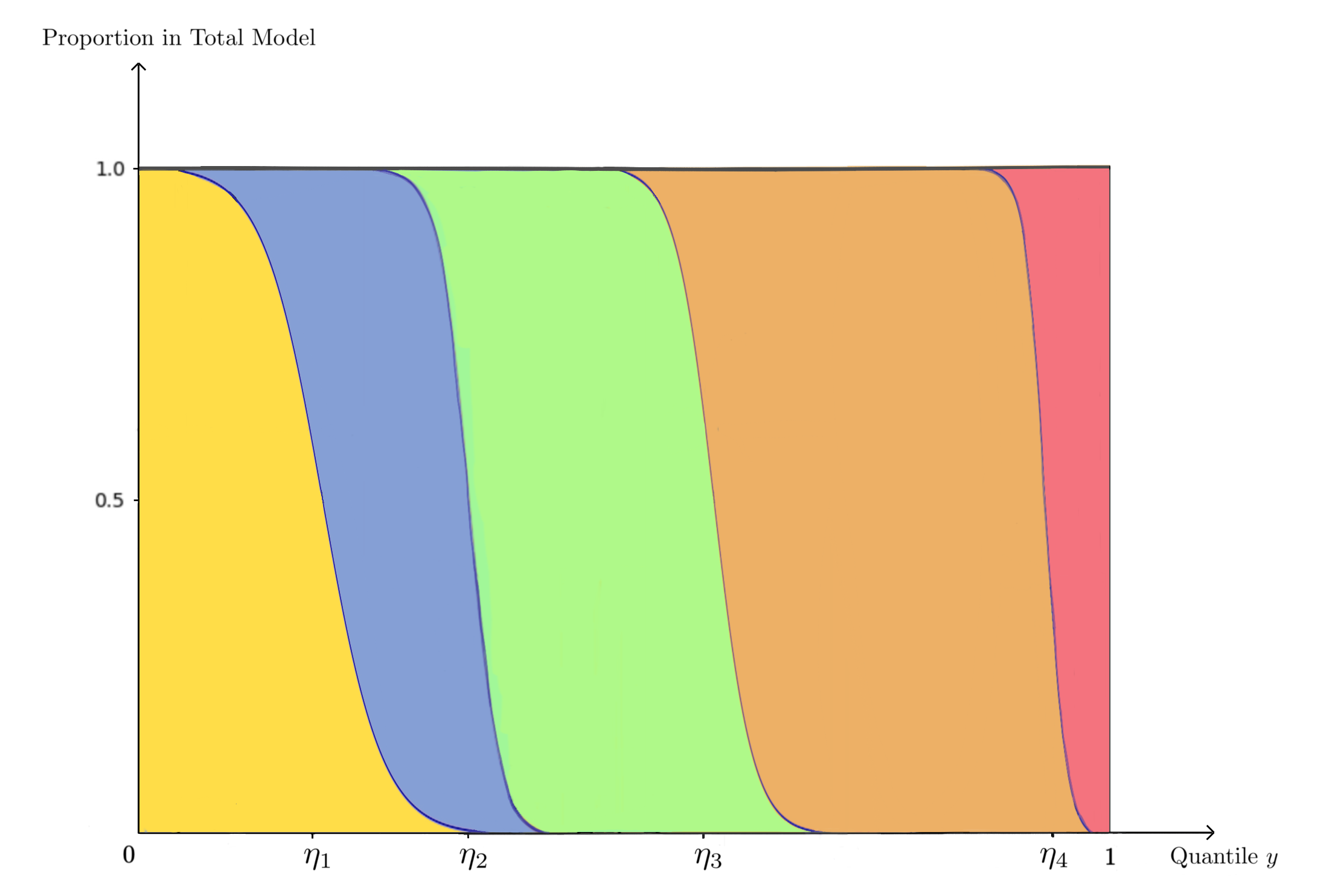}
    \caption{Visual Illustration of the Smoothed $K$-segment Model ($K=5$).}
    \label{modplot}
\end{figure}

Figure~\ref{modplot} visually depicts the influence of the smoothing method on each submodel's proportion within the smoothed $K$-segment model, specifically illustrated for $K=5$. The x-axis is the quantile $y$ with the segment separation at $\eta_k$ for $k=1,...,4$. The y-axis is the portion of submodel $S_k$ within the overall $K$-segment model for each $y$ value, and each submodel is distinguished by color. 

We propose two specific smoothing methods later in this section: quantile-based smoothing and score-based smoothing.
quantile-based smoothing employs the two submodels nearest to a segmentation boundary, connecting them via a modified version of the sigmoid function to reestimate valuations; score-based smoothing calculates contributions from all submodels, adjusting their weights based on the distance from the sales value.
%Each of these methods offers unique benefits in refining the model.

\subsubsection*{Quantile-based Smoothing}
%\red{We should rethink the name ``Quantile-based smoothing'' and ``score-based smoothing''.}

The basic idea of quantile-based smoothing is to use a shifted and scaled sigmoid function to smooth out the boundary between two adjacent submodels. More specifically, for any $k=1,...,K-1$, we define 
\[g_k(y) = \sigma\left(-\frac{10}{\gamma_k-\eta_k+\lambda_k}(y-\gamma_k)-5\right), \quad \text{for } y \in [\eta_k-\lambda_k,\gamma_k], \]
where $\sigma(z)=1/(1 + e^{-z})$ is the sigmoid function, $\lambda_k$ and $\gamma_k$ are the hyper-parameters of the function $g_k(y)$. $g_k(y)$ is a continuously differentiable S-shaped function monotonically decreasing from $\sigma(5)$ to $\sigma(-5)$ on its domain $[\eta_k-\lambda_k,\gamma_k]$.\footnote{\(\sigma(-5) < 0.01\) and \(\sigma(5) > 0.99\). Thus, for practical applications, one can effectively consider the range of the function $g_k(y)$ as $[0,1]$.}

% The domain of $g_k(y)$ is $[\eta_k-\lambda_k,\gamma_k]$, and the range of $g_k(y)$ belongs to $[0,1]$. 
% Roughly speaking, $g_k(y)$ is a continuously differentiable function from near $1$ to near $0$.
% we introduce the function $g_k(y)$ which acts on value intervals around the segment boundary $\eta_k$ for $k=1,\ldots,K-1$.
% In particular, we define $g_k(y)$ as a continuously differentiable function, adapted from the sigmoid function $\sigma(z)=1/(1 + e^{-z})$.  
% We define smoothing parameters $\lambda_k$ and $\gamma_k$ for $g_k$, which outline the interval where \( g_k(y) \) operates, specifically marking the left endpoint at \( \eta_k - \lambda_k \) and the right endpoint at \( \gamma_k \).
% The parameters adhere to the relation $\gamma_{k-1}+\lambda_k\leq\eta_k\leq\gamma_k$ for $ k=1,\ldots,K-1$. For notation consistency, we set $\gamma_0=0$. 
% The function \(g_k\) capitalizes on an S-shaped curve similar to the sigmoid function \(\sigma\), which effectively maps any input to a range between 0 and 1. 
% , which facilitates a smooth transition between adjacent submodels.
% For a sales price within the range $[\eta_k-\lambda_k,\gamma_k]$, particularly near $\eta_k$, the smoothing method will combine the influences of both the $k$-th and $(k+1)$-th submodels, using the smoothing function $g_k$. As the sales price moves away from the segmentation boundary, its valuation is driven by the submodel of the segment $y$ falls within.
The \textbf{Quantile-based Smoothing $K$-segment Model} is defined as
\[M_{\text{q}}(y)=
\begin{cases}
\begin{aligned}
& S_k(y), && \text{if } \gamma_{k-1} \leq y < \eta_k-\lambda_k ,\\
& g_k(y)\cdot S_k(y)+\big(1-g_k(y)\big)\cdot S_{k+1}(y), && \text{if } \eta_k-\lambda_k\leq y< \eta_{k},\\
& S_k(y), && \text{if } \eta_{K-1}\leq y\leq 1.
\end{aligned}
\end{cases}
\]
Roughly speaking, $M_{\text{q}}(y)$ aims to use the $k$-th submodel when the quantile $y$ falls into $[\eta_{k-1}, \eta_k]$ and, at the same time, provide a smoothing separation around the boundary between two adjunct segments. The width of the smoothing separation is given by the parameters $\lambda_k$ and $\gamma_k$, which are hyper-parameters to be tuned in the submodels. Notice that the quantile-based smoothing method introduces $2K$ hyper-parameters, which can potentially lead to overfitting with fine-tuning. To address this issue, we introduce score-based smoothing.

 % Employing the two submodels adjacent to a segmentation boundary yields two significant advantages. First, it enhances the interpretability of the model’s specific valuations by maintaining a clear and intuitive structure. Second, it incorporates a set of tunable parameters $\{(\lambda_k,\gamma_k)\}_{k=1}^{K-1}$ that can be tailored to balance the impact of adjacent submodels on the final estimtae.  

% However, over-relying on fine-tuning parameters for specific local data may make the quantile-based smoothing approach more complex, potentially leading to overfitting, and neglect of information in distant segments. This could introduce systematically skewed valuations.
% To address these concerns, we propose another smoothing technique that incorporates a broader spectrum of submodels for valuation, streamlines parameter requirements, and ensures a uniform and interpretable model expression across all segments.

\subsubsection*{Score-based Smoothing}

Instead of considering only one submodel, the idea of \textbf{Score-based Smoothing} utilizes the outputs of all submodels, and assigns a weight to each of the submodels, i.e.,
\begin{equation}\label{weightingsmoothing}
    M_{\text{WS}}(y) = \sum_{k=1}^K w_k(y) S_k(y),
\end{equation}
where $w_k(y)$ is the weight of the $k$-th submodel for a quantile $y$. For each submodel \(S_k\), the weight \( w_k(y) = s_k(y)/\sum_{j=1}^K s_j(y) \), is derived by normalizing a score function, \( s_k(y) \). We here present two approaches to define the score:
%\red{The name ``monotonic'' needs a second thought.}
\begin{itemize}
    \item \textbf{Midpoint score.} The midpoint score for the $k$-th segment is defined as
    \begin{equation}\label{nomono}
    s_k(y) = e^{-\frac{\mu}{\eta_k-\eta_{k-1}}\cdot\big|y-\frac{\eta_k+\eta_{k-1}}{2}\big|}\ ,
    \end{equation}
    where $\mu$ is a hyper-parameter of the model.
    \item \textbf{Distance score.} The distance score for the $k$-th segment is defined as
    \begin{equation}\label{mono}
    s_k(y) = e^{-\frac{\mu}{\eta_k-\eta_{k-1}}\cdot\operatorname{dist}\big(y,[\eta_{k-1},\eta_k]\big)},
    \end{equation}
    where $\mu$ is a hyper-parameter of the model, and \(\operatorname{dist}(y, [\eta_{k-1},\eta_k]) = \inf_{a \in [\eta_{k-1},\eta_k]} \|y - a\|\).
\end{itemize}

Both scores decrease exponentially with the distance to the $k$-th segment $[\eta_{k-1}, \eta_k]$. The difference between them is that the midpoint score utilizes the distance to the midpoint of the segment, i.e., $\frac{1}{2}(\eta_{k-1}+\eta_k)$, while the distance score utilizes the distance to the set $[\eta_{k-1}, \eta_k]$. Intuitively, if $y$ is inside a segment, the corresponding score is high and, if it is far away, the corresponding score is low. The hyper-parameter $\mu$ controls how sensitive to the distance to a segment the score is.

\begin{figure}[H]
\centering  
\subfigure[$\textit{M}_{\textit{m--s},\textit{5}}$]{
\includegraphics[width=0.45\textwidth]{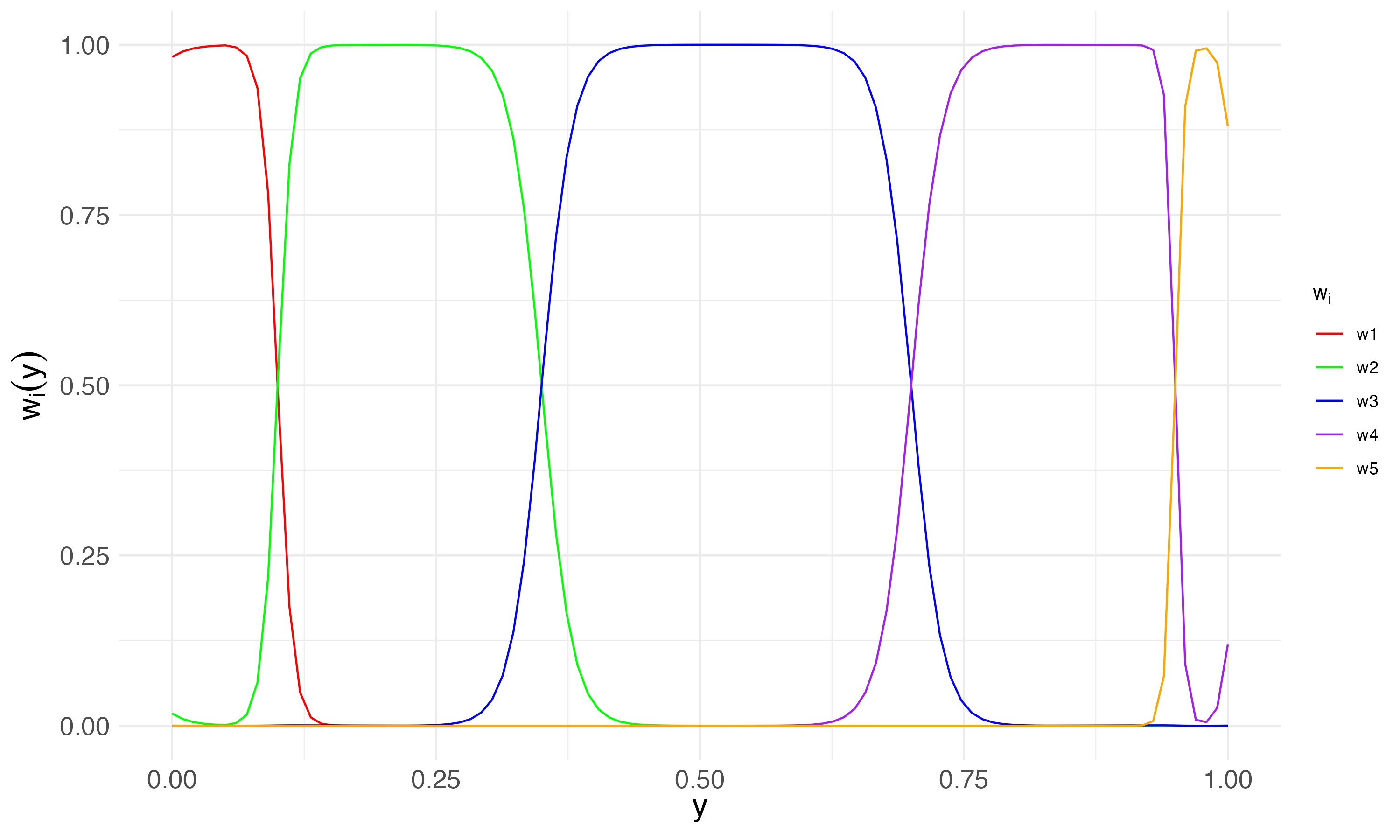}}
\subfigure[$\textit{M}_{\textit{d--s},\textit{5}}$]{
\includegraphics[width=0.45\textwidth]{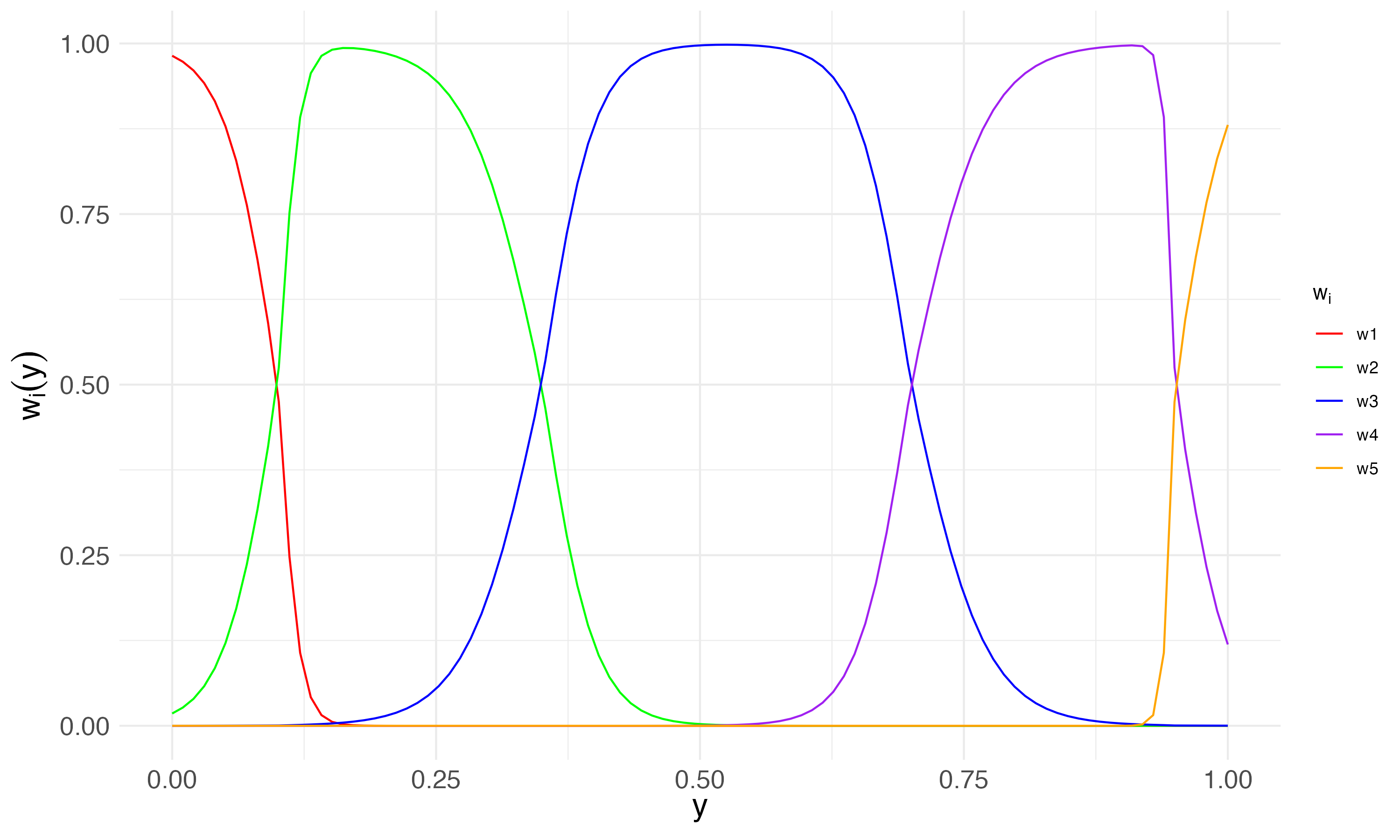}}
\caption{Weighting function $w_i(y)$ versus quantile $y$ (when $K=5,\eta_1=0.1,\eta_2=0.35,\eta_3=0.7,\eta_4=0.95,\mu=10$).}
\label{k5}
\end{figure}

 As an illustration, with \(K=5\), Figure~\ref{k5} depicts the relationship between \(y\) and the respective weights \(w_k\) for \(k=1,\ldots,5\), as defined by score function $s_k(y)$ in Equations~(\ref{nomono}) and (\ref{mono}). Panel (a) of Figure~\ref{k5} shows a 'non-monotonic' pattern in the corresponding weight of the midpoint score, where the weight $w_i$ sometimes decreases within the $i$-th segment and increases when moving away from this segment. Distance score shown in panel (b), effectively resolves this problem and leads to a ``monotonic" pattern.

% \begin{figure}[H]
% \centering  
% \subfigure[$w_1(y)$ against  $y$]{
% \includegraphics[width=0.3\textwidth]{plot7/curvenonmono_1.png}}
% \subfigure[$w_2(y)$ against  $y$]{
% \includegraphics[width=0.3\textwidth]{plot7/curvenonmono_2.png}}
% \subfigure[$w_3(y)$ against  $y$]{
% \includegraphics[width=0.3\textwidth]{plot7/curvemono_3.png}}
% \subfigure[$w_4(y)$ against  $y$]{
% \includegraphics[width=0.3\textwidth]{plot7/curvenonmono_4.png}}
% \subfigure[$w_5(y)$ against  $y$]{
% \includegraphics[width=0.3\textwidth]{plot7/curvenonmono_5.png}}
% \caption{Weighting Function of $M_{\text{m-s},5}$ when $\eta_1=0.1,\eta_2=0.35,\eta_3=0.7,\eta_4=0.95,\mu=10$.}
% \label{k5nm}
% \end{figure}

% \begin{figure}[H]
% \centering  
% \subfigure[$w_1(y)$ against  $y$]{
% \includegraphics[width=0.3\textwidth]{plot7/curvemono_1.png}}
% \subfigure[$w_2(y)$ against  $y$]{
% \includegraphics[width=0.3\textwidth]{plot7/curvemono_2.png}}
% \subfigure[$w_3(y)$ against  $y$]{
% \includegraphics[width=0.3\textwidth]{plot7/curvemono_3.png}}
% \subfigure[$w_4(y)$ against  $y$]{
% \includegraphics[width=0.3\textwidth]{plot7/curvemono_4.png}}
% \subfigure[$w_5(y)$ against  $y$]{
% \includegraphics[width=0.3\textwidth]{plot7/curvemono_5.png}}
% \caption{Weighting Function of $M_{\text{d-s},5}$ when $\eta_1=0.1,\eta_2=0.35,\eta_3=0.7,\eta_4=0.95,\mu=10$.}
% \label{k5}
% \end{figure}

We end this section by discussing two major advantages of the $K$-segment model over a single-segment model.
\begin{enumerate}
    \item Flexibility: 
The $K$-segment model divides data into multiple segments, enabling adaptation to complex, nonlinear relationships within each segment. Instead of imposing a single structure that assumes uniformity throughout the entire dataset, the $K$-segment model recognizes that different segments of data can exhibit unique behaviors. This provides flexibility in understanding unique property value characteristics, while a one-size-fits-all model might lead to oversimplification.
\item 
Increased Accuracy: By training individual submodels on more homogeneous data within each segment, the $K$-segment model mitigates the impact of skewed data. The specific tailoring of each submodel to its respective segment enables a more precise understanding of property values within that range. This method leads to a more accurate assessment across varied property values, surpassing a single-segment model that applies a uniform approach to all properties, regardless of value differences.
\end{enumerate}

\section{Performance Evaluation}\label{eva}

In this section, we conduct numerical experiments to evaluate the performance of the $K$-segment model. Our findings indicate that it effectively mitigates regressivity, leading to enhanced accuracy and fairness in property taxation. Section~\ref{data} details the dataset employed for our experiments. In Section~\ref{par}, we specify the parameters of our models. 
In Section~\ref{se1}, we numerically show an improvement of predicted $R^2$ obtained by the $K$-segment model relative to the original model. 
In Section~\ref{plot}, we illustrate the regressivity observed in the original model's assessment, and how the $K$-segment model can help to address this issue.
Section~\ref{se2} presents a comparison of accuracy and fairness between the various models.
In Section~\ref{se3}, we analyze the trade-off between accuracy and fairness in our models. 

% Here, we amalgamate insights from the prior phases, contrasting the accuracy and fairness attributes of the $K$-segment model under various smoothing techniques and parameters, leading to the identification of the most advantageous model configuration.

\subsection{Dataset and Submodel Training} \label{data}

% \textbf{Dataset:} 
Our experiments use two primary data sources from CCAO \cite{CCAO2022}: \texttt{training\_data} and \texttt{assessment\_data}. The \texttt{training\_data} comprises residential sales from nine years before the assessment year (2013--2021) and serves as the foundation for training and evaluating our model. It consists of 410,584 rows, each representing an observation, and $167$ columns,\footnote{Both \texttt{training\_data} and \texttt{assessment\_data}'s columns include basic information and features of properties. Due to recording discrepancies, the number of columns differs between the two datasets. However, all essential features and variables used during model training and assessment are encompassed within datasets.} each representing a feature used in property assessments or related recorded data.
And the \texttt{assessment\_data} has more samples, comprising 1,099,111 rows and $165$ columns. 
%\footnote{The sales prices within the \texttt{assessment\_data} span from \$37,771.95 to \$2,931,788.} 
Specifically, we construct the following datasets following~\cite{CCAO2022}:
\begin{enumerate}
    \item \textbf{Training set}: It comprises the initial 90\% of the \texttt{training\_data}. This set is for training the submodels, which can be further divided into:
         \begin{itemize}
            \item Training Subset: Primarily used for model fitting.
            \item Validation Subset: Used for optimizing the hyper-parameters of the LightGBM model. 
         \end{itemize}
    This subdivision follows the rolling-origin resampling technique as described in \cite{hyndman2018forecasting}, expanding the training subset size within a stable time window and designating the subsequent period as the validation set, which ensures robust hyper-parameter selection. See \cite{CCAO2022,hyndman2018forecasting} for more details.
    \item \textbf{Test set}: It contains the last 10\% of \texttt{training\_data}. This set evaluates the performance of the refined model on unseen data. Once the desired performance is achieved, the model is retrained using the entirety of the \texttt{training\_data}.
    \item \textbf{Assessment set}: This is the \texttt{assessment\_data}, which includes all properties for assessment and is where the fully trained model is applied.
\end{enumerate}

%\subsection{Submodel training}
For each segment, we use the original model, as described in Section~\ref{modelfor} and further detailed in \cite{CCAO2022}, to develop submodels from the \texttt{training\_data}.\footnote{To simplify the workflow, and reduce potential overfitting, we chose not to use the post-modeling phase in either the original model or the submodels of the $K$-segment model.}
These submodels leverage LightGBM for training, with its parameters tuned through Bayesian hyper-parameter optimization as discussed in \cite{kuhn2019feature,CCAO2022}.
The configurations for the parameters of the $K$-segment model can be found in Section~\ref{par}.
Then, we use the trained model to assess property values in the \texttt{assessment\_data}.

% We follow the methodology outlined in \cite{CCAO2022} to develop submodels using \texttt{training\_data}. From \texttt{training\_data}, the Train set is derived from the initial 90\% of sales data, while the Test set consists of the most recent 10\%.
% The Train set is then split further into a training subset and a validation subset, purposed for model fitting and hyper-parameter optimization of the LightGBM model. This subdivision follows the rolling-origin resampling technique as described by \cite{hyndman2018forecasting}, expanding the training set size within a stable time window and designating the subsequent period as the validation set, which ensures robust hyper-parameter selection. 
% Once satisfactory performance measures are achieved on the Test set, the model is retrained using the entirety of the Train set, including the validation subset. Our subsequent work and contributions will build on this foundation.
% Subsequently, the Assessment set (\texttt{assessment\_data}), which covers all residential properties for valuation, is analyzed using our final model to produce the property value predictions.

% In summary, we will focus on the following sets in our analysis:
% \begin{itemize}
%     \item Train set: Used for training submodels and choosing the hyper-parameter combination to ensure optimal accuracy.
%     \item Test set: Used to assess the performance of the trained and fine-tuned model on unseen data.
%     \item Assessment set: The dataset of properties to be assessed, where the fully trained model is applied.
% \end{itemize}

\subsection{Parameter Configuration}\label{par}

In this section, we 
specify the parameters used in the \(K\)-segment model, which includes the number of segments $K$, the segmentation parameters $\{\eta_k\}_{k=1}^K$, the smoothing parameters \(\{\lambda_k,\gamma_k\}_{k=1}^K\) for the quantile-based method, and the smoothing parameter \(\mu\) for the score-based smoothing method.

The number of segments, denoted by $K$, strongly impacts the performance of the $K$-segment model. 
A $K$-segment model with a small \(K\) value may behave similarly to the original model, negating the benefits of segmentation.
While a larger \(K\) enhances the model's expressivity, it also introduces overfitting and heightened computational demands, as each segment requires individual cross-validation for the optimal tuning of LightGBM hyper-parameters.
We here opted for $K=3,5$ in our experiments, which we believe provides a balance between these conflicting factors.

Table \ref{tab:table2} details the model parameters used in this numerical study. After selecting $K$, we set the segmentation parameters $\eta$ correspondingly as shown in Table \ref{tab:table1}. Then we test the four smoothing methods described in Section \ref{model}. For quantile-based smoothing, we also specify the parameters $\gamma$ and $\lambda$. For score-based smoothing, we set $\mu=10$. These parameters are tuned to achieve the best $R^2$ for the test set.

\begin{table}[H]
\label{tab:table1}
\centering
\begin{tabular}{c|c|c}
               & Segmentation Parameters & Smoothing Methods                \\ \hline
Original Model & -                       & -                                 \\ \hline
\multirow{4}{*}{3-segment Model} & \multirow{4}{*}{\begin{tabular}[c]{@{}c@{}c@{}}$(\eta_1, \eta_2) = (0.1, 0.9)$\\ $(\lambda_1,\lambda_2)=(0.1, 0.1)$\\$(\gamma_1,\gamma_2)=(0.2, 1)$  \end{tabular}}                    & Unsmoothed                          \\ \cline{3-3} 
               &                         & Quantile-based Smoothing        \\ \cline{3-3} 
               &                         & Midpoint Score-based Smoothing \\ \cline{3-3} 
               &                         & Distance Score-based Smoothing     \\ \hline
\multirow{4}{*}{5-segment Model}   & \multirow{4}{*}{\begin{tabular}[c]{@{}c@{}c@{}}$(\eta_1, \eta_2, \eta_3, \eta_4) = (0.2, 0.35, 0.7, 0.9)$\\$(\lambda_1,\lambda_2,\lambda_3,\lambda_4)=(0.15, 0.03, 0.1, 0.1)$\\$(\gamma_1,\gamma_2,\gamma_3,\gamma_4)=(0.3, 0.5, 0.73, 1)$  \end{tabular}}    & Unsmoothed                          \\ \cline{3-3} 
               &                         & Quantile-based Smoothing        \\ \cline{3-3} 
               &                         & Midpoint Score-based Smoothing \\ \cline{3-3} 
               &                         & Distance Score-based Smoothing    
\end{tabular}
\caption{Parameters and Smoothing Methods for the $K$-segment Model ($K = 3, K = 5$).}
\end{table}

\subsection{Accuracy: \texorpdfstring{$R^2$}{R square} of Models}\label{se1}

Table~\ref{r2} compares the $R^2$ of the $K$-segment and original models. A few observations are in order:
\begin{enumerate}
    \item The $K$-segment model consistently improves accuracy, as reflected in the $R^2$ values, especially for the assessment set.\footnote{The differing $R^2$ value ranges between the unseen test set and assessment set is attributed to the partial availability of 2022 sales price data at the time of our research.}
    \item The $K$-segment model, for both $K=3$ and $K=5$, show similar $R^2$ improvements over the original model, suggesting comparable performance between them.
    \item Introducing smoothing methods further boosts the accuracy of the $K$-segment model. While the extent of improvement in $R^2$ varies depending on the smoothing method and the fixed $K$ value, smoothed models outperform their unsmoothed counterparts in all experiments.
\end{enumerate}

\begin{table}[H]
\centering
\setlength{\tabcolsep}{7mm}{
\begin{tabular}{c|ccc}
\rowcolor[HTML]{FFFFC7} 
\cellcolor[HTML]{FFFFC7} &
  \multicolumn{3}{c}{\cellcolor[HTML]{FFFFC7}$R^2$} \\ \cline{2-4} 
\multirow{-2}{*}{\cellcolor[HTML]{FFFFC7}Model} &
  \multicolumn{1}{c|}{Train} &
  \multicolumn{1}{c|}{Test} &
  Assessment \\ \hline
  \text{$M_{unsm,3}$} &
  \multicolumn{1}{c|}{0.9900958} &
  \multicolumn{1}{c|}{0.8559155} &
  0.7622086\\
   \rowcolor[HTML]{67FD9A}   \text{$M_{unsm,5}$} &
  \multicolumn{1}{c|}{\cellcolor[HTML]{67FD9A}0.992362} &
  \multicolumn{1}{c|}{\cellcolor[HTML]{67FD9A}0.8486896} &
 0.7692979\\
   \text{$M_{q,3}$} &
  \multicolumn{1}{c|}{0.9870508} &
  \multicolumn{1}{c|}{0.8583131} &
  0.7811874\\
 \rowcolor[HTML]{67FD9A}   \text{$M_{q,5}$} &
  \multicolumn{1}{c|}{\cellcolor[HTML]{67FD9A}0.9841513} &
  \multicolumn{1}{c|}{\cellcolor[HTML]{67FD9A}0.8651624} &
  0.7750398\\
  \text{$M_{m-s,3}$} &
  \multicolumn{1}{c|}{0.9803078} &
  \multicolumn{1}{c|}{0.8388794} &
  0.7975622\\
 \rowcolor[HTML]{67FD9A}   \text{$M_{m-s,5}$} &
  \multicolumn{1}{c|}{\cellcolor[HTML]{67FD9A}0.9988296} &
  \multicolumn{1}{c|}{\cellcolor[HTML]{67FD9A}0.8533657} &
  0.7797632\\
 \text{$M_{d-s,3}$} &
  \multicolumn{1}{c|}{0.9806621} &
  \multicolumn{1}{c|}{0.837238} &
  0.7845766\\
  % \text{$M_{d-s,3(B)}$} &
  % \multicolumn{1}{c|}{0.9871261} &
  % \multicolumn{1}{c|}{0.8476007} &
  % 0.7886364\\
 \rowcolor[HTML]{67FD9A}   \text{$M_{d-s,5}$} &
  \multicolumn{1}{c|}{\cellcolor[HTML]{67FD9A}0.9975327} &
  \multicolumn{1}{c|}{\cellcolor[HTML]{67FD9A}0.8553439} &
  0.7892524\\
\rowcolor[HTML]{ECF4FF} 
Original & \multicolumn{1}{c|}{\cellcolor[HTML]{ECF4FF}0.9863093} & \multicolumn{1}{c|}{\cellcolor[HTML]{ECF4FF}0.8290886} & 0.708173
\end{tabular}}
\caption{$R^2$ Performance of Original and $K$-segment Models for Different Datasets.\protect\footnotemark}
\label{r2}
\end{table}
\footnotetext{In this and the following tables, entries for the original model are highlighted in blue, and those for the $K$-segment model with $K=5$ are marked in green for ease of differentiation.}

% Initially, we discuss the unsmoothed $K$-segment model, establishing a base performance prior to the introduction of smoothing techniques. The original model, while fitting the Train set closely with an $R^2$ of 0.986, shows limited generalizability to the Test and Assessment sets. The unsmoothed $K$-segment model imprvoes both accuracy and generalizability.
% When $K=3$, the model raises the $R^2$ from 0.829 to 0.856 (a 3.3\% increase) on the Test set, and from 0.708 to 0.762 (a 7.6\% increase) on the Assessment set. For $K=5$, the $R^2$ value on the Test set increases from 0.829 to 0.849 (a 2.4\% increase), and on the Assessment set it escalates from 0.708 to 0.769 (an 8.6\% increase).

% The improvement in $R^2$ values for $K=3$ and $K=5$ is similar, although the $K$-segment model with $K=5$ provides tighter predictions per sales price as discussed in Section~\ref{3.2}.

% Next, we turn to the impact of smoothing methods on the $K$-segment model. For $K=3$, the best-performing smoothing method improves the $R^2$ from 0.762 (unsmoothed) to 0.798 (midpoint), marking an increase of 4.7\%. For $K=5$, the best-performing smoothing methods raises the $R^2$ from 0.769 (unsmoothed) to 0.789 (Distance), making a growth of 2.6\%. While the exact $R^2$ improvement differs between different smoothing methods and parameters for a given $K$ value, all methods yield an improvement compared to the unsmoothed model.

\subsection{Regressivity: Visual Representations}\label{plot}

\begin{figure}[H]
\centering  
\subfigure[Original Model]{
\includegraphics[width=0.18\textwidth]{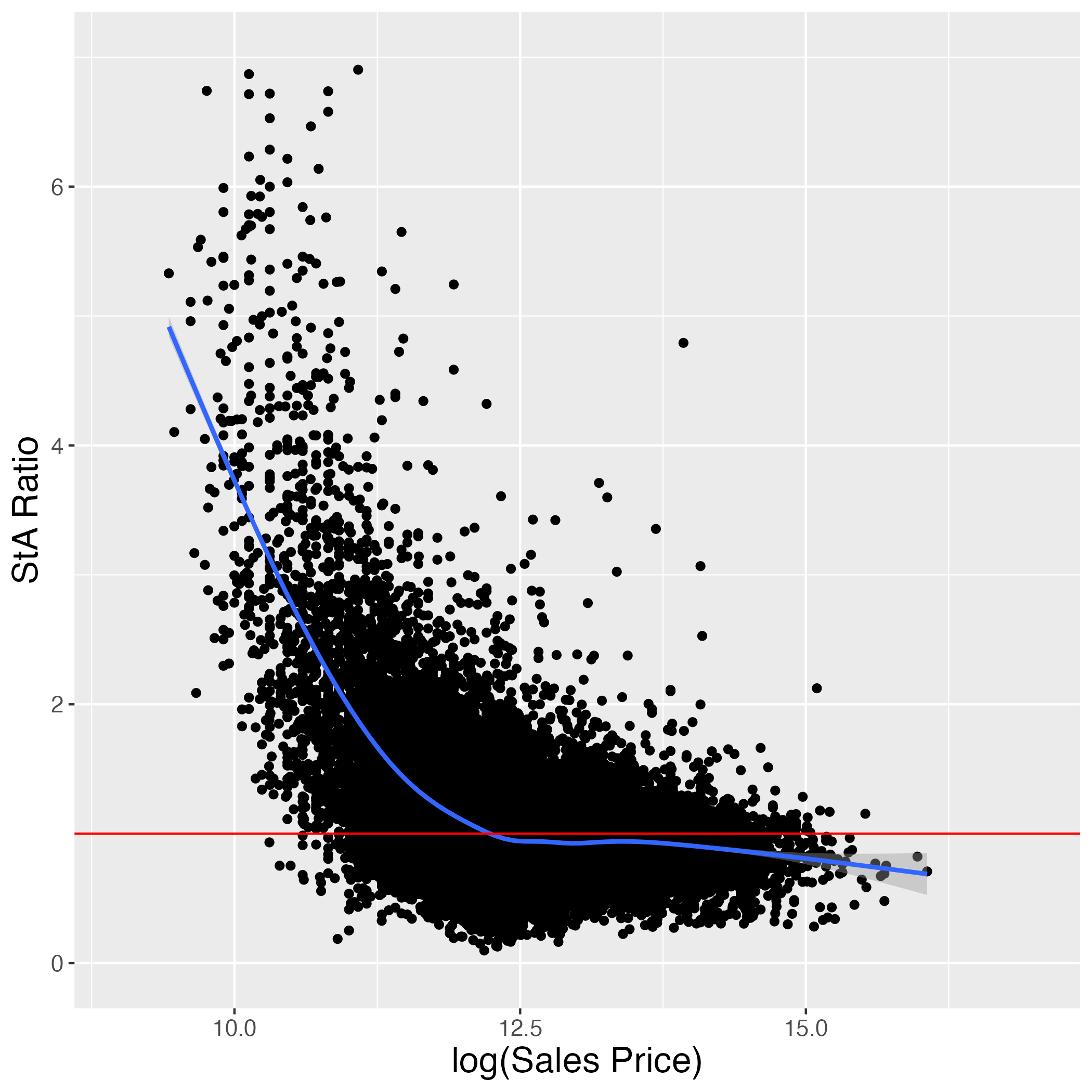}}
\subfigure[$\textit{M}_{\textit{d--s},\textit{3}}$]{
\includegraphics[width=0.18\textwidth]{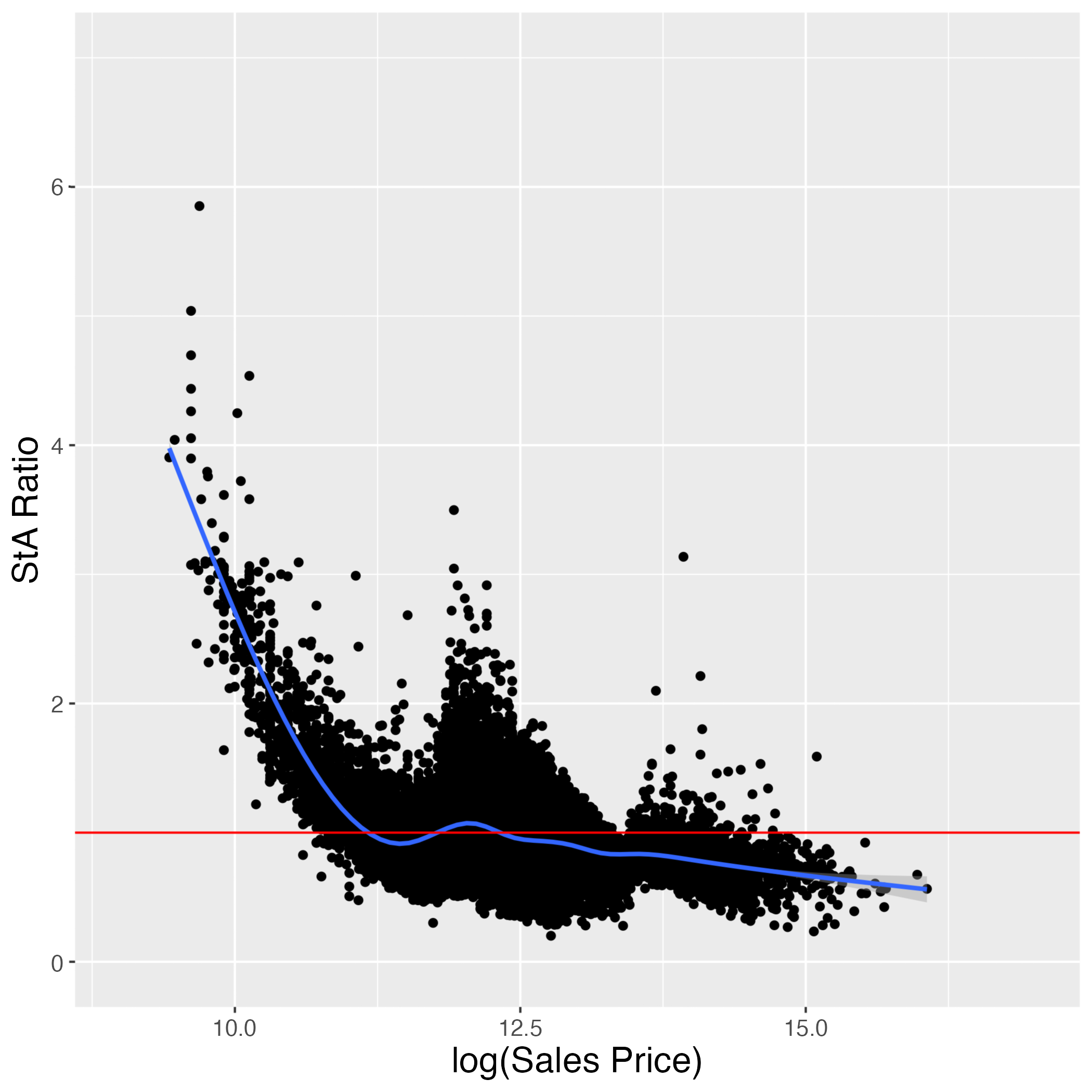}}
\subfigure[$\textit{M}_{\textit{q},\textit{5}}$]{
\includegraphics[width=0.18\textwidth]{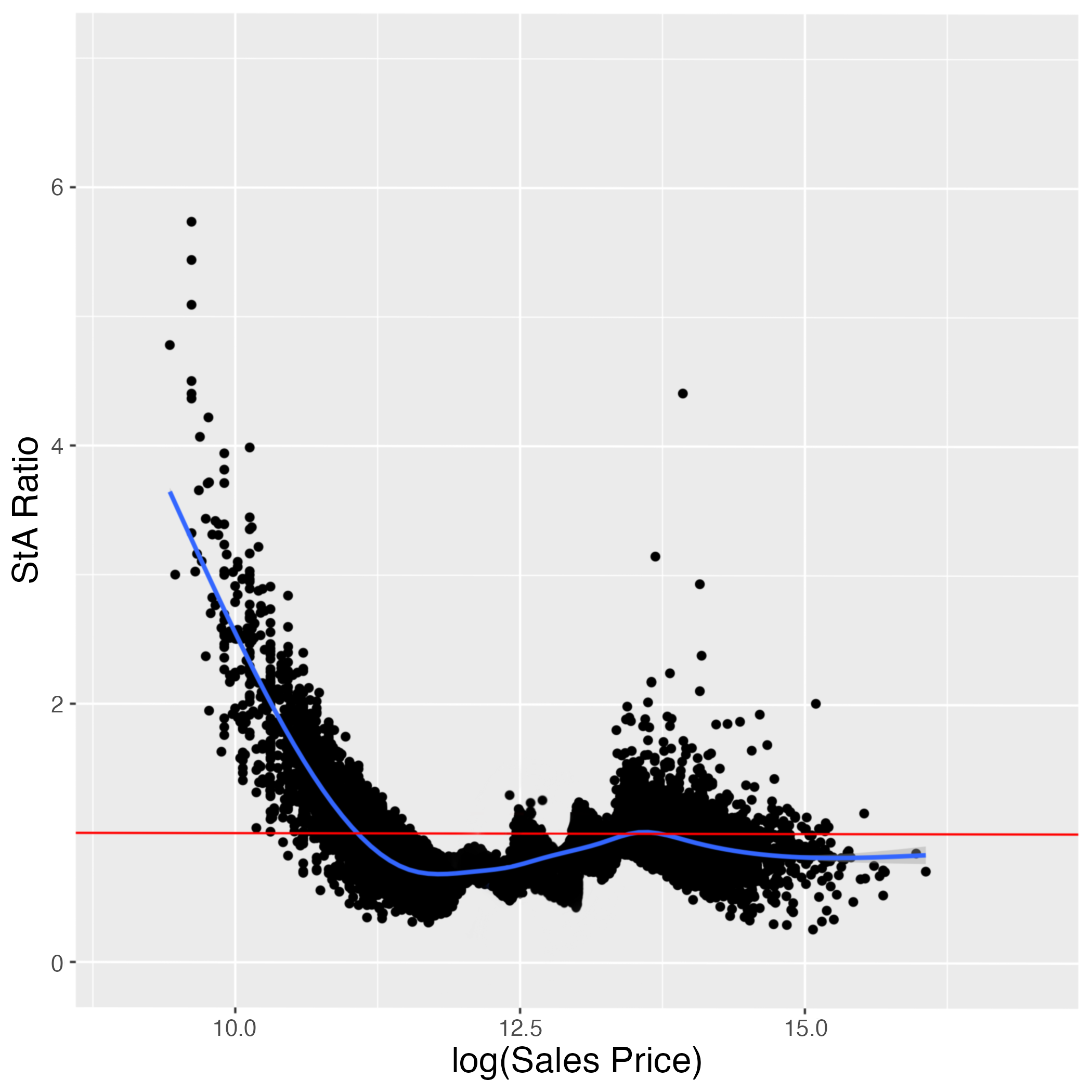}}
\subfigure[$\textit{M}_{\textit{m--s},\textit{5}}$]{
\includegraphics[width=0.18\textwidth]{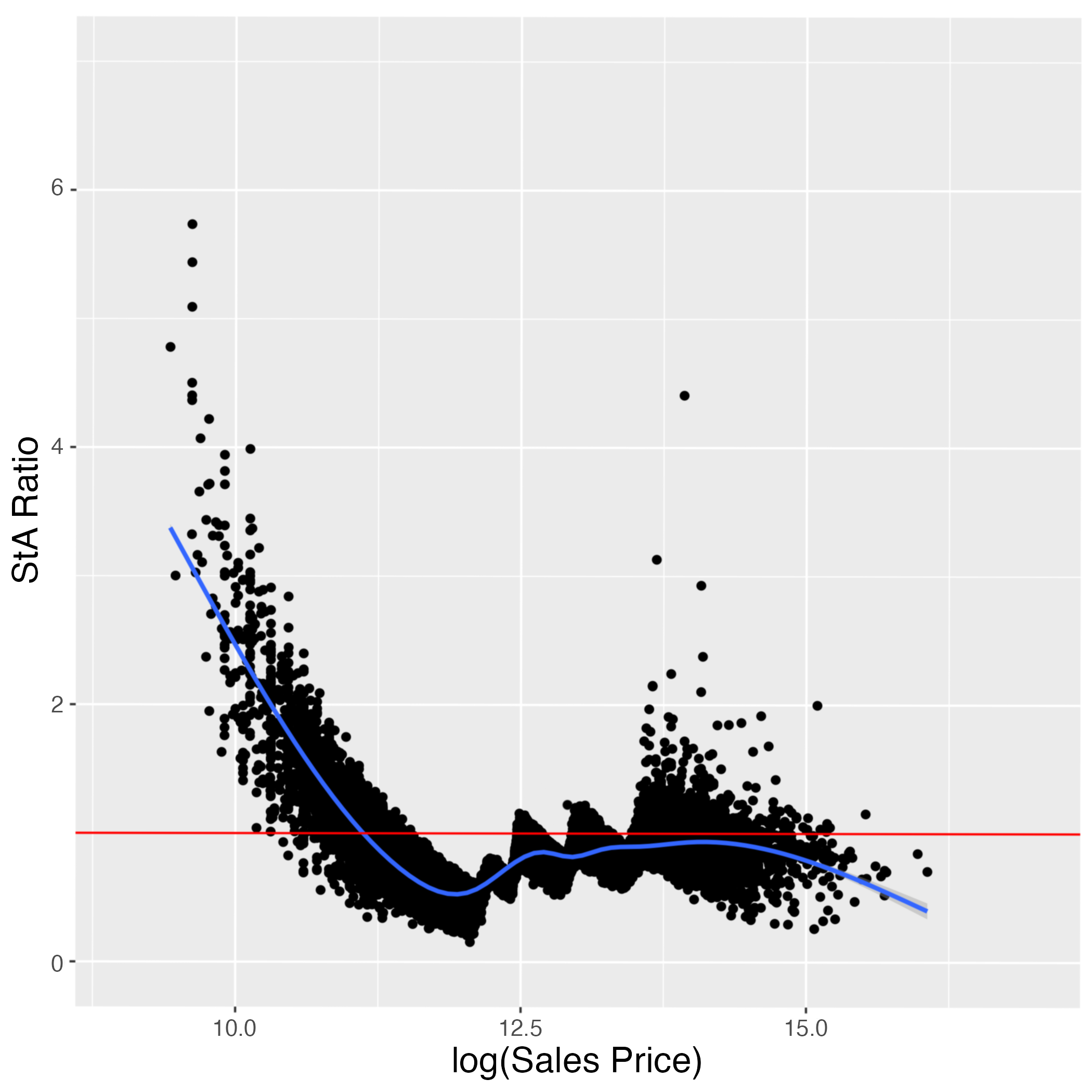}}
\subfigure[$\textit{M}_{\textit{d--s},\textit{5}}$]{
\includegraphics[width=0.18\textwidth]{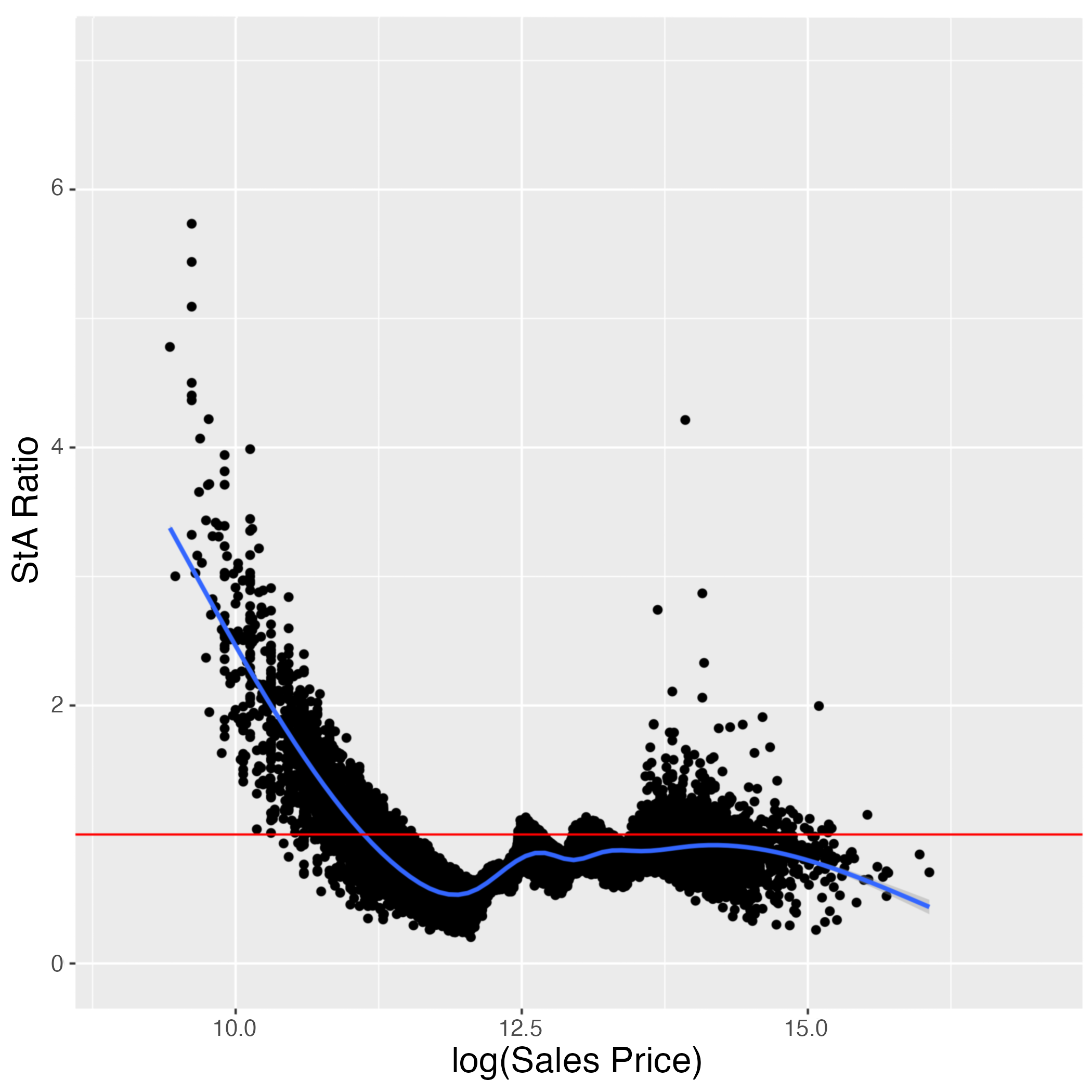}}
\caption{StA Ratio ($r$) versus Logarithm of Sales Price ($\log(x)$).}
\label{ratio2}
\end{figure}

In this section, we compare the regressivity of the $K$-segment and original models. Figure~\ref{ratio2} presents the relationship between the sales prices and the StA ratio of different models, using available data from properties sold in 2021. In each of the panels of Figure~\ref{ratio2}, the x-axis is the logarithm of the sales price $(\log(x))$, and the y-axis is the StA ratio. Each dot is a sample in the corresponding dataset.  The red and blue lines indicate the baseline StA ratio (equal to 1) and the trend line obtained from a Generalized Additive Model (GAM) using the data, respectively.\footnote{GAMs are used for smoothing in this context due to their ability to model complex, nonlinear relationships. They serve as the default method in \texttt{ggplot2} and provide a useful reference for visualizing trends. However, the choice of smoothing method can influence the interpretation of the plot.} Following are a few observations:\footnote{We concentrate on the price range that the majority of properties fall within, excluding outliers, and analyze sales with logarithm values contained within $[9,16]$.}

% \begin{figure}[H]
% \centering  
% \subfigure[Original Model]{
% \label{Fig.new2.sub.1}
% \includegraphics[width=0.18\textwidth]{Redefined-segment-Model/plots4/ori_logratio_assess.png}}
% \subfigure[$M_{\text{d-s},3}$]{
% \label{Fig.new2.sub.2}
% \includegraphics[width=0.18\textwidth]{Redefined-segment-Model/plot0/ori_logratio_assess.png}}
% \subfigure[$M_{\text{q},5}$]{
% \label{Fig.new2.sub.2}
% \includegraphics[width=0.18\textwidth]{Redefined-segment-Model/plot5/5seg/ori_logratio_assess.png}}
% \subfigure[$M_{\text{m-s},5}$ ]{
% \label{Fig.new2.sub.2}
% \includegraphics[width=0.18\textwidth]{Redefined-segment-Model/plot6/5seg/ori_logratio_assess.png}}
% \subfigure[$M_{\text{d-s},5}$]{
% \includegraphics[width=0.18\textwidth]{plot8/ori_logratio_assess.png}}
% \caption{StA Ratio ($r$) versus Logarithm of Sales Price ($\log(x)$) in the Assessment Set.}
% \label{ratio3}
% \end{figure}

\begin{itemize}
    \item The original model shows clear regressivity, as can be seen from the negative correlation between the StA ratio and the sales price, shown in panel (a) of Figure~\ref{ratio2}.
    % In the Assessment set, due to potential market changes from the Covid-19 pandemic, a larger proportion of valued homes have StA ratios above 1.
    % \item The $K$-segment models notably mitigate regressivity, in particular on the train and the test dataset.
    % \item As $K$ increases in the $K$-segment model, the StA ratio adjustment becomes more noticeable, with $K=5$ showing a significant mitigation in regressivity.
    \item The $K$-segment model increases data point concentration and, as $K$ increases, this concentration intensifies.
    \item Smoothing models, given a fixed $K$, efficiently smooth segments, with only subtle distinctions between them.
\end{itemize}

% \subsubsection{Valuation Comparisons: \texorpdfstring{$K$}{K}-Segment Model vs. Original Model}\label{3.3}
\begin{figure}[t]
\centering 
\subfigure[$\textit{M}_{\textit{q},\textit{3}}$]{
\includegraphics[width=0.45\textwidth]{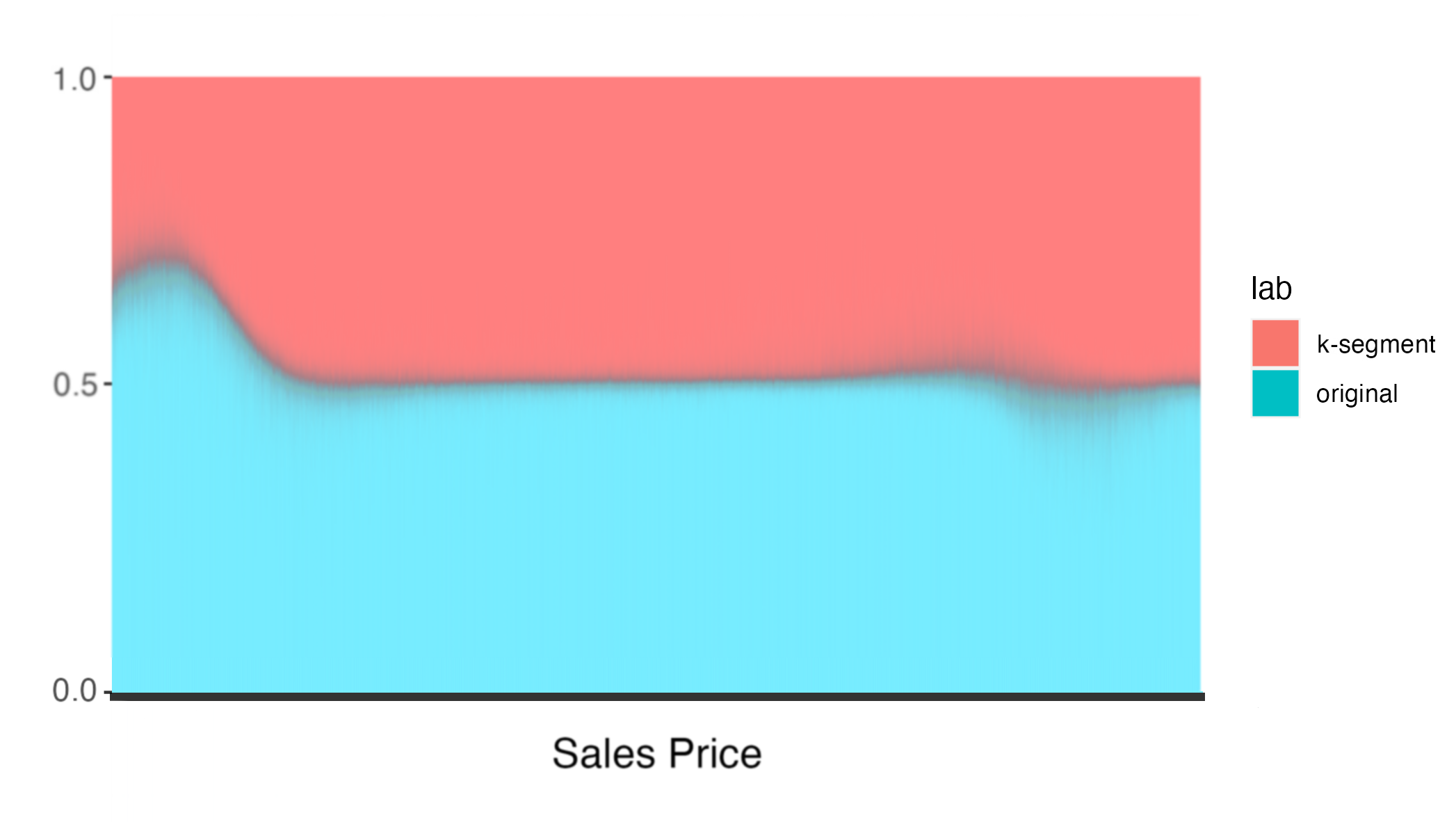}}
\subfigure[$\textit{M}_{\textit{q},\textit{5}}$]{
\includegraphics[width=0.45\textwidth]{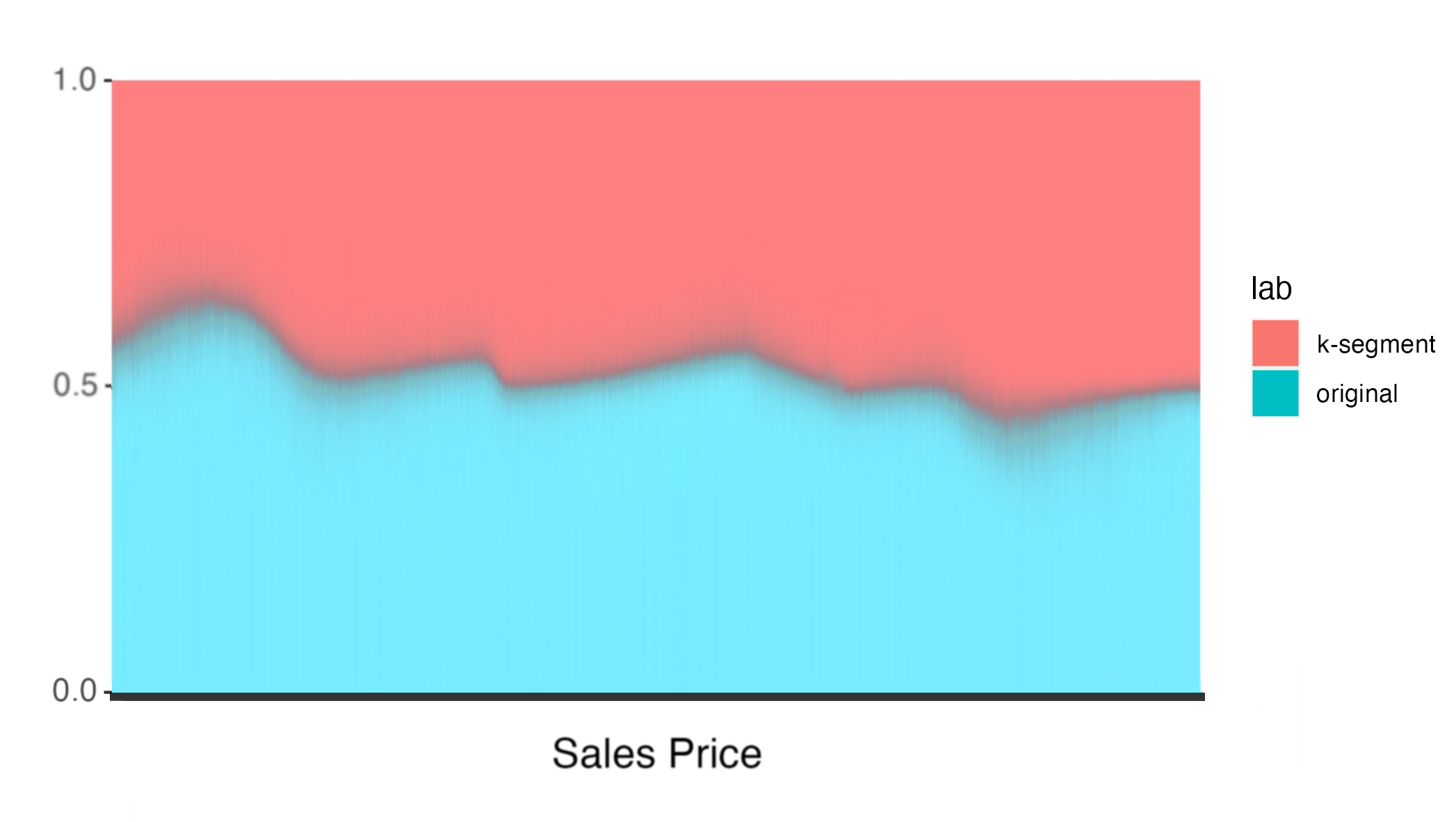}}
\subfigure[$\textit{M}_{\textit{m--s},\textit{3}}$]{
\includegraphics[width=0.45\textwidth]{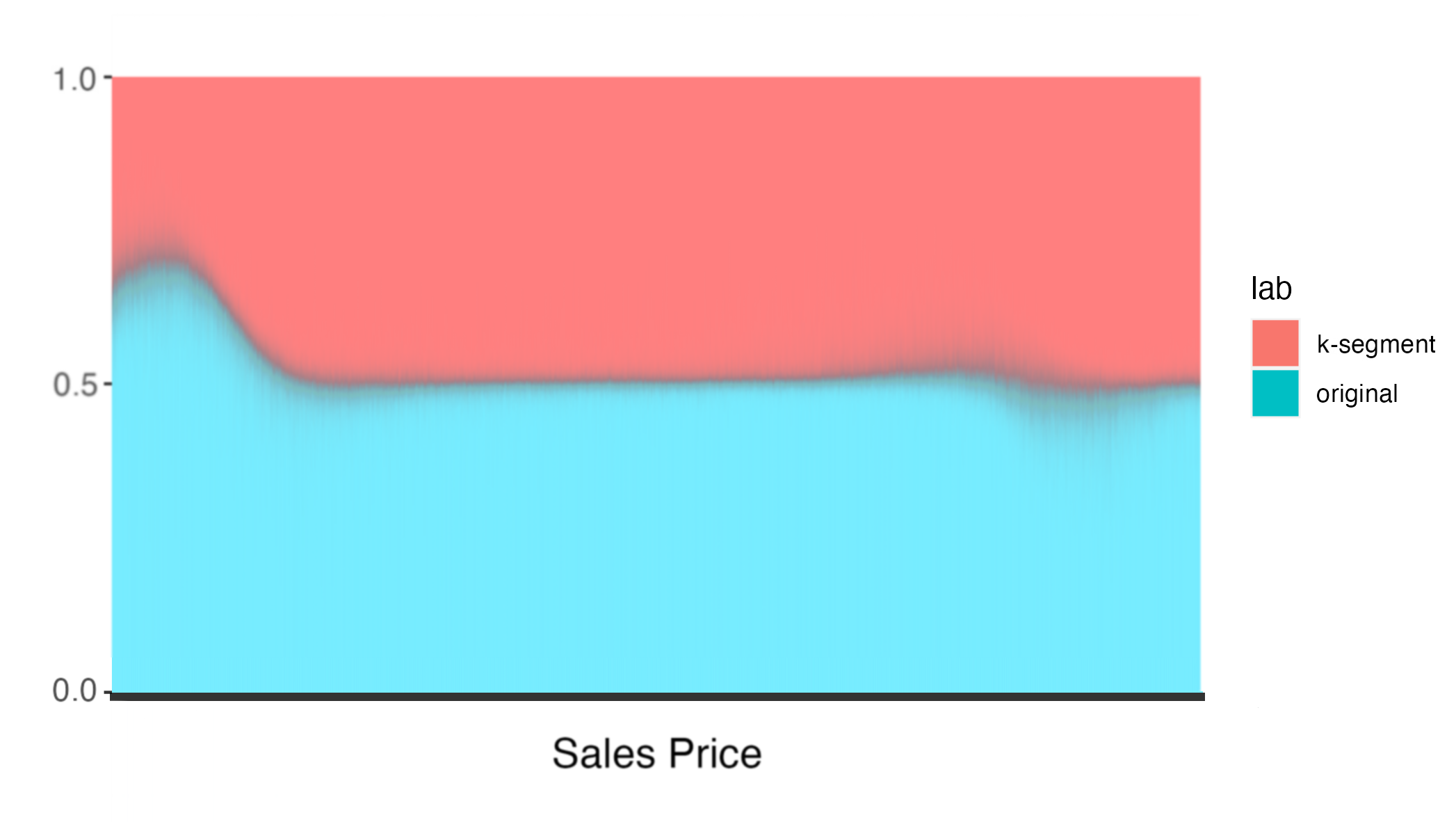}}
\subfigure[$\textit{M}_{\textit{m--s},\textit{5}}$]{
\includegraphics[width=0.45\textwidth]{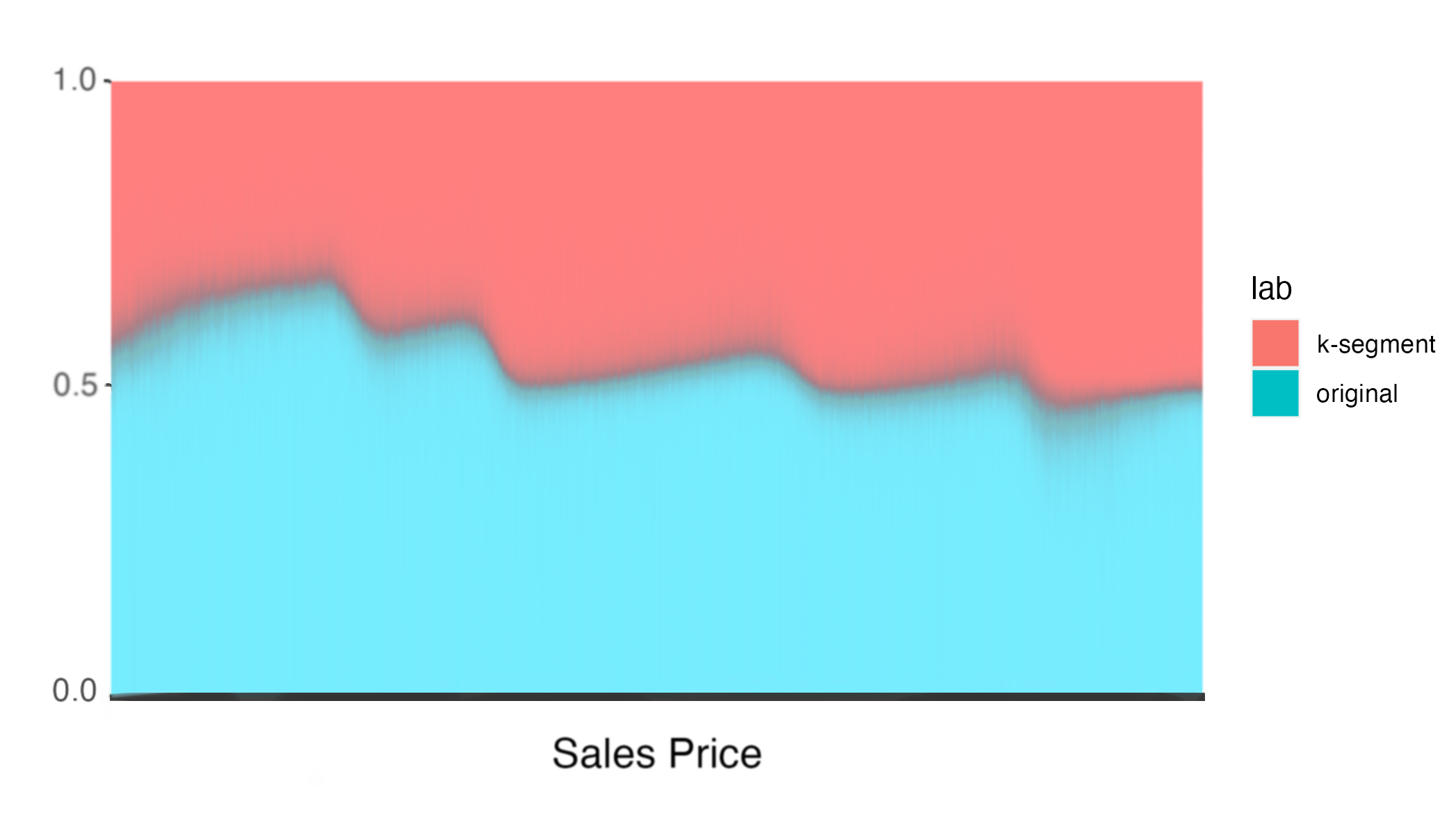}}
\subfigure[$\textit{M}_{\textit{d--s},\textit{3}}$]{
\includegraphics[width=0.45\textwidth]{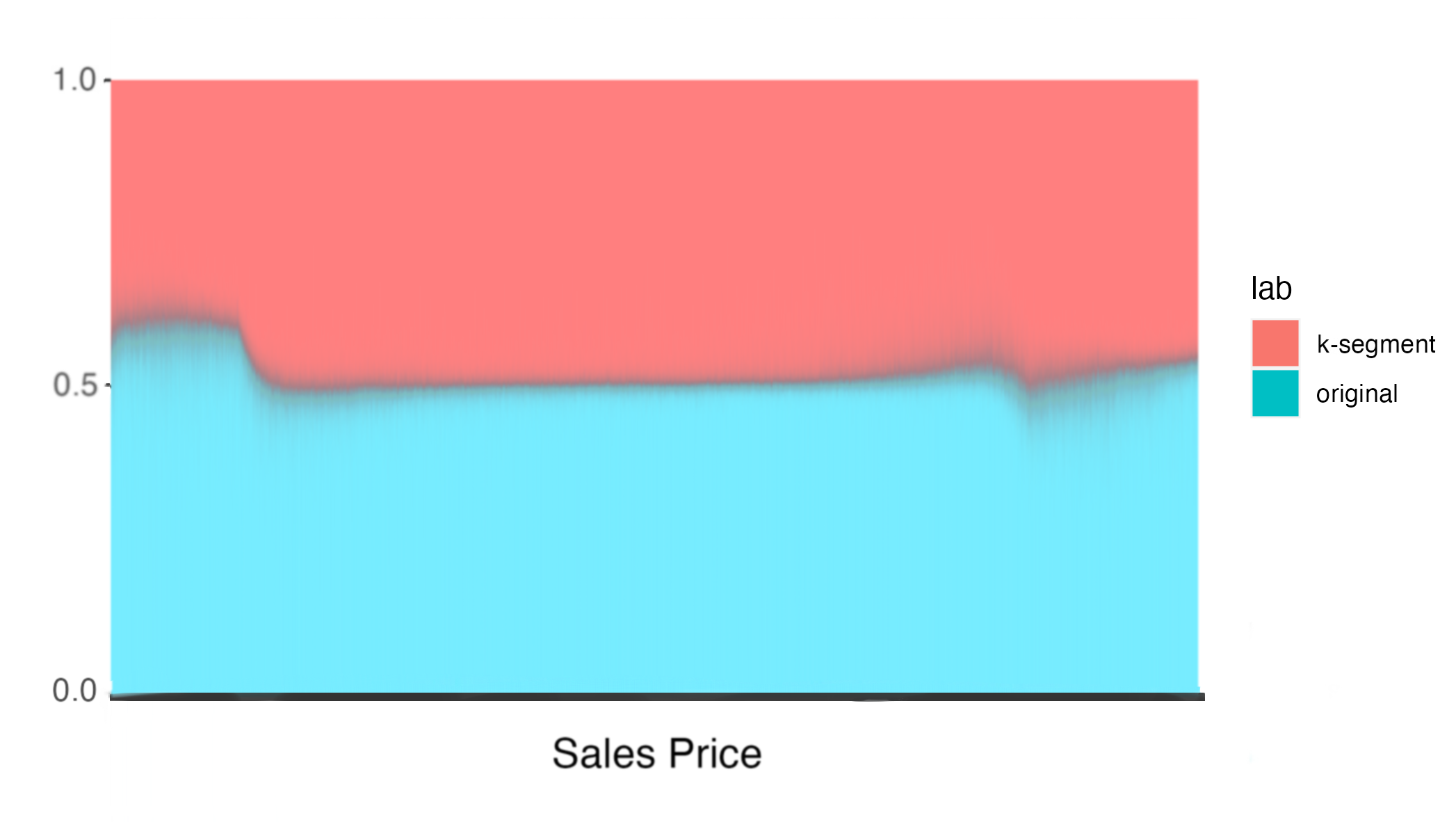}}
\subfigure[$\textit{M}_{\textit{d--s},\textit{5}}$]{
\includegraphics[width=0.45\textwidth]{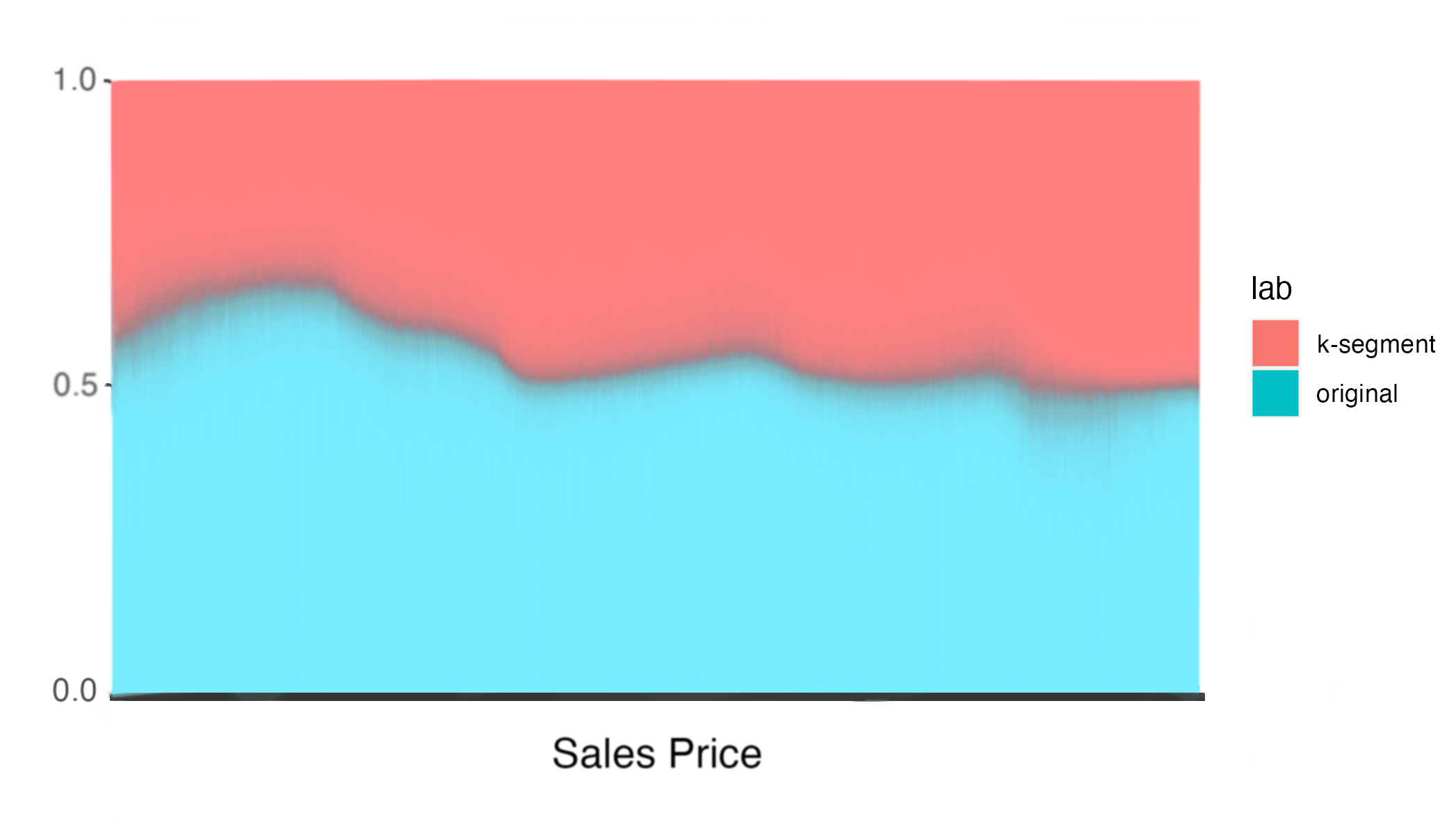}}
\caption{Comparison of the $K$-segment and Original Models by Sales Price.}
\label{per}
\end{figure}

Figure~\ref{per} compares assessments of the $K$-segment and original models by sales price levels, using available data from properties sold in 2021.
The x-axis represents sales prices, which are categorized into levels and arranged in ascending order.
The y-axis shows a comparison of the evaluations of the $K$-segment model (with different parameters in (a)--(f)) and the original model.
Specifically, for each sales price, 
\[\frac{\text{Length of red bar}}{\text{Length of blue bar}}=\frac{\text{Assessment price of } K\text{-segment model}}{\text{Assessment price of original model}}.\]
% The x-axis depicts sales prices in a relative format by classifying each sales price into distinct ordered levels.
% The y-axis contrasts the predictions of the original model (blue) and the $K$-segment model (red), with their combined height standardized to 1.
% More specificaly, Let's denote the assessment value of the original model as \( A_{\text{original}} \) and the assessment value of the $K$-segment model as \( A_K\) for a given sales price level.
% The blue segment's length is given by $A_{\text{original}}/\big(A_{\text{original}} + A_{K}\big) $. The red segment's length is given by $A_{K}/\big(A_{\text{original}} + A_{K}\big) $. A blue section above the 0.5 mark indicates that the $K$-segment model predicts a lower assessment for that price level compared to the original model.
We can observe from Figure \ref{per} that

\begin{itemize}
    \item Compared to the original model, the $K$-segment model tends to produce lower estimates for lower-valued properties and higher estimates for higher-valued properties; thus it can mitigate regressivity. 
    \item The score-based smoothing with Distance score method can smooth out the rapid change around the boundary between submodels (see panels (e) and (f) of Figure~\ref{per}).
    % The smoothness, relative to other K-segment models, is evident in Figure~\ref{Valuation Index} (e) and Figure~\ref{per} (f), where the Distance score-based smoothing method is applied.
\end{itemize}

%\red{Rewrite this paragraph.} Figure~\ref{per} presents a dual-colored bar plot, to visually contrast the assessments of the $K$-segment and original models for varying sales prices in the Assessment set. The x-axis represents sales prices. The y-axis shows the relative averaged valuations of the two models (red for the $K$-segment model and blue for the original model), displayed as proportions of their combined total. This highlights how the $K$-segment model's assessments vary from the original model across different sales prices.
%\footnote{Categorizing sales prices aids in emphasizing distinct property value ranges, especially at both spectrum ends. It improves readability by preventing plot skewness that might result from a higher density of lower-valued homes. This method ensures a balanced, interpretable visualization without being overshadowed by any specific price range.}

To summarize, Figures~\ref{ratio2} and~\ref{per} show that the original model exhibits regressivity, and that the $K$-segment model reduces regressivity by producing lower estimates for lower-valued properties and higher estimates for higher-valued properties.

% \begin{enumerate}
%     \item 
%     \item 
% \end{enumerate}

\subsection{Fairness: Group Fairness and Deviation-weighted Fairness}\label{se2}

In the section, we apply the fairness measures outlined in Section~\ref{fair} to the original and $K$-segment models for properties sold in 2021. Using the original model as a baseline for comparison, we illustrate the superiority of the fairness measure of the $K$-segment model.
%parameterized as detailed in Section~\ref{par}.

\subsubsection{Group Fairness}\label{sectiongrp}

Group fairness, as characterized in Definition~\ref{grp}, quantifies the discrepancies in the StA ratios between distinct groups. These groups are differentiated based on sales price percentiles, encapsulating varied property value ranges. 
A larger disparity, indicating that properties across different price percentiles are not assessed uniformly, yields a lower $F_{grp}$ value, which denotes a lack of fairness. Conversely, a smaller disparity, resulting in a higher $F_{grp}$ value, signifies that the assessments are more uniform and therefore fairer across groups.

In our experimental setup, we implement a two-step configuration for group fairness. First, we set the number of groups, $n$, to be either 2 or 3.\footnote{We choose $n=2$ and $n=3$ for group fairness configurations in our experiment due to their practical effectiveness and interpretability. When $n=2$, we compare each pair of StA ratios between the ``low" and ``high" groups, focusing on the extremes of the property valuation distribution. For $n=3$, we assess the difference in StA ratios between the ``low", ``medium", and ``high" groups, providing a more comprehensive evaluation. Both values balance the granularity of the analysis and the complexity of interpreting the results.} Second, we segment the groups by dividing the data into equally sized percentiles. The group fairness measure ($F_{grp}$) employed tends to produce relatively small absolute values (see Appendix~\ref{app3} for details). Here, we introduce the concept of relative unfairness to rescale group fairness measure using the unfairness score of the original model:  %Specifically, we set $m_t=\lfloor m/n\rfloor$ for $t=1,\ldots, n-1$. 
%This methodology ensures that each group encompasses a balanced range of property values, allowing a comprehensive evaluation of fairness across all price tiers.

% In Figure~\ref{grpfair}, the x-axis represents $F_{grp}$ with $n=2$, and the y-axis is $F_{grp}$ with $n=3$. The following are observations.

% Different $K$ show substantial impact on group fairness, with larger $K$ values generally resulting in higher $F_{grp}$ scores. 
%Meanwhile, the nearly linear relationship between $n=2$ and $n=3$ suggests that the relative performance of different models remains consistent across varying $n$ values. This visualization also aids in selecting the optimal smoothing method to enhance group fairness. Both midpoint and Distance weighted smoothing methods in the $K$-segment models with $K=5$ achieve superior group fairness scores for various values $n$.

\begin{definition}\label{uf}
[Relative Unfairness]
The relative unfairness of group fairness, denoted by $RU_{grp}$, is obtained by comparing the $F_{grp}$ value of a given model $M$ with that of the original model. Specifically, it is calculated as:
         \[RU_{grp}=\frac{F_{grp}\text{ of model }M}{F_{grp} \text{ of original model}}.
         \]
\end{definition}

% \begin{figure}[H]
% \centering  
% \subfigure[$F_{grp}$ Across Parameters: $n=2$ vs. $n=3$.]{
% \label{grpfair}
% \includegraphics[width=0.45\textwidth]{plot8/groupfairness.png}}
% \subfigure[$F_{dev}$ Across Parameters: $\alpha=2$ vs. $\alpha=5$.)]{
% \label{dynfair}
% \includegraphics[width=0.45\textwidth]{plot8/dynfairness.png}}
% \caption{Comparison of Fairness Measures Across Parameters.}
% \end{figure}

The \( RU_{grp} \) measure facilitates fairness comparisons between models. Set against a baseline where the original model has an \( RU_{grp} \) of 1, values greater than 1 suggest unfairness, while those below 1 suggest fairness. The smaller the value, the fairer the model. Table~\ref{table2.5} compares models based on their \( RU_{grp} \) scores. As we can see,

\begin{table}[H]
\centering
\begin{tabular}{>{\columncolor[HTML]{FFFFC7}}c l|c
>{\columncolor[HTML]{67FD9A}}c c
>{\columncolor[HTML]{67FD9A}}c c
>{\columncolor[HTML]{67FD9A}}c cc
>{\columncolor[HTML]{67FD9A}}c >{\columncolor[HTML]{ECF4FF}}c }
\multicolumn{2}{c|}{\cellcolor[HTML]{FFFFC7}Model}  &
  $M_{unsm,3}$ &
  {\color[HTML]{333333}$M_{unsm,5}$} &
  $M_{q,3}$ &
  $M_{q,5}$ &
  $M_{m-s,3}$ 
  % &
  % $M_{m-s,5}$ &
  % $M_{d-s,3}$ &
  % $M_{d-s,3(B)}$ &
  % $M_{d-s,5} $&
  % $\text{Ori}$ 
  \\ \hline
 \multicolumn{1}{c|}{\cellcolor[HTML]{FFFFC7}}&
  $n=2$ &
  0.83 &
  {\color[HTML]{333333}0.19} &
  0.70 &
  0.15 &
  0.78 
  % &
  % 0.14 &
  % 0.86 &
  % 0.74 &
  % 0.13 &
  % 1 
\\ \cline{2-7} 
\multicolumn{1}{c|}{\multirow{-2}{*}{\cellcolor[HTML]{FFFFC7}Relative Unfairness of $F_{grp}$}}&
  $n=3$ &
  0.81 &
  {\color[HTML]{333333}0.15} &
  0.68 &
  0.17 &
  0.76 
  % &
  % 0.10 &
  % 0.85 &
  % 0.69 &
  % 0.11 &
  % 1
\end{tabular}
\\
\vspace{0.2cm}
\begin{tabular}{>{\columncolor[HTML]{FFFFC7}}c l|
>{\columncolor[HTML]{67FD9A}}c c
>{\columncolor[HTML]{67FD9A}}c >{\columncolor[HTML]{ECF4FF}}c }
\multicolumn{2}{c|}{\cellcolor[HTML]{FFFFC7}Model}  &
  {\color[HTML]{333333}$M_{m-s,5}$} &
  $M_{d-s,3}$ &
  % $M_{d-s,3(B)}$ &
  $M_{d-s,5} $&
  $\text{Original}$ 
  \\ \hline
 \multicolumn{1}{c|}{\cellcolor[HTML]{FFFFC7}} &
  $n=2$ &
  {\color[HTML]{333333}0.14} &
  0.86 &
 % 0.74 &
  0.13 &
  1 
\\ \cline{2-6} 
\multicolumn{1}{c|}{\multirow{-2}{*}{\cellcolor[HTML]{FFFFC7}Relative Unfairness of $F_{grp}$}} &
  $n=3$ &
  {\color[HTML]{333333}0.10} &
  0.85 &
%  0.69 &
  0.11 &
  1
\end{tabular}
\caption{Relative Unfairness of Group Fairness in the Original and $K$-segment Models.}
\label{table2.5}
\end{table}

\begin{enumerate}
\item The $K$-segment model positively influences group fairness. While the improvement is modest for $K=3$, it becomes  more significant for $K=5$;
\item Smoothing methods can enhance group fairness. However, the optimal smoothing method can vary depending on the chosen $K$ value.
% \item Larger $n$ values lead to increased group fairness scores. However, $RU_{grp}$ remains largely unaffected for $n=2,3$\footnote{Since the $(\cdot)^+$ operation would yield zero when the StA ratios of lower-valued homes is smaller or equal to that of higher-valued homes. If no regressivity existed in the property assessment, there shouold be a marked discrepancy from the theoretical $9/2$ times relationship between the scores for $n=3$ and $n=2$. However, our data exhibit only a minor deviation from the $9/2$ ratio. This suggests that in most cases, the StA ratios for lower-valued homes (lower groups) are indeed higher than those of the higher-valued homes (higher groups). This minor deviation illustrates the regressive pattern in our data source and underscores the necessity for addressing this trend to achieve fairer property assessments.}.
\end{enumerate}

\subsubsection{Deviation-weighted Fairness}\label{sectiondev}

Deviation-weighted fairness, as introduced in Definition~\ref{dyn}, uses weighting functions to emphasize the overestimations of lower-valued properties and underestimations of higher-valued properties. 
The optimal state of this fairness measure is achieved when the StA ratio is consistently 1 for all properties, indicating accurate and uniform valuations. $F_{dev}$ decreases with increasing deviations from the optimal StA ratio.
% Consequently, a model with a higher $F_{dev}$ score represents better fairness due to fewer deviations from the optimal StA ratio.

In the computation of $F_{dev}$, we select an $\alpha$ value from the set $\{0,1,2,5\}$. The selection is based on the degree of variation these values induce in $e^{-\alpha \Tilde{y}}$ for $\Tilde{y} \in [0,1]$. With a larger $\alpha$, the weighting function $\Tilde{w}_1(\Tilde{y})$ decreases more rapidly as $\Tilde{y}$ increases. 

% Figure~\ref{dynfair} visualizes how deviation-weighted fairness can be improved using $K$-segment models, and it also provides guidance for choosing optimal smoothing methods. 
%The x-axis of Figure~\ref{dynfair} indicates $F_{dev}$ with $\alpha=2$, whereas the y-axis indicates $F_{dev}$ with $\alpha=5$. The following are some observations.

% It demonstrates that models with larger $K$ values typically have higher $F_{dev}$ scores, suggesting their potential to improve fairness.
%It also reveals a positive correlation between the outcomes for $\alpha=2$ and $\alpha=5$, indicating that a model's performance is relatively consistent for different $\alpha$ values when $\alpha$ is sufficiently large\footnote{However, it is worth noting that this consistency may not hold when $\alpha$ is smaller, as seen in the case of $\alpha=0$. When $\alpha$ is small, the value of $e^{-\alpha y}$ is closer to 1 for all $y \in [0,1]$, which means that the weights $w_1(y)$ and $w_2(y)$ do not vary as significantly across different property values. This results in the deviation-weighted fairness $F_{dev}$ being less sensitive to detect these regressivity issues. Therefore, different models, even those that may differ significantly in their degree of regressivity, can end up with similar $F_{dev}$ values.}.
% The highest $F_{dev}$ scores are notably achieved by $K$-segment models with $K=5$, utilizing Quantile-based smoothing methods.

The absolute values of $F_{dev}$ scores are generally quite large (see Appendix~\ref{app3;2}). To better compare these scores across different models, we introduce a relative unfairness measure for deviation-weighted fairness, denoted by $RU_{dev}$, that is formulated as follows:
 \[RU_{dev}=\frac{F_{dev}\text{ of model }M}{F_{dev} \text{ of original model}},\] by comparing the $F_{dev}$ values of a specific model $M$ and the original model. The smaller the $RU_{dev}$, the fairer the model.
Table~\ref{table3.5} compares the relative unfairness, $RU_{dev}$, of different models, which shows that

\begin{table}[H]
\centering
\begin{tabular}{>{\columncolor[HTML]{FFFFC7}}c l|c
>{\columncolor[HTML]{67FD9A}}c c
>{\columncolor[HTML]{67FD9A}}c c
>{\columncolor[HTML]{67FD9A}}c cc
>{\columncolor[HTML]{67FD9A}}c >{\columncolor[HTML]{ECF4FF}}c }
\multicolumn{2}{c|}{\cellcolor[HTML]{FFFFC7}Model}  &
  $M_{unsm,3}$ &
  {\color[HTML]{333333}$M_{unsm,5}$} &
  $M_{q,3}$ &
  $M_{q,5}$ &
  $M_{m-s,3}$ 
  % &
  % $M_{m-s,5}$ &
  % $M_{d-s,3}$ &
  % $M_{d-s,3(B)}$ &
  % $M_{d-s,5} $&
  % $\text{Ori}$ 
  \\ \hline
\multicolumn{1}{c|}{\cellcolor[HTML]{FFFFC7}} &
  $\alpha=0$ &
1.07 &
  {\color[HTML]{333333}1.12} &
  0.95 &
  0.71 &
  0.99 
\\ \cline{2-7} 
\multicolumn{1}{c|}{\cellcolor[HTML]{FFFFC7}}&
  $\alpha=1$ &
  0.97 &
  {\color[HTML]{333333}0.91} &
  0.86 &
  0.55 &
  0.91 \\ \cline{2-7} 
\multicolumn{1}{c|}{\cellcolor[HTML]{FFFFC7}}&
  $\alpha=2$ &
  0.90 &
  {\color[HTML]{333333}0.74} &
  0.78 &
  0.43 &
  0.84 
   \\ \cline{2-7} 
\multicolumn{1}{c|}{\multirow{-4}{*}{\cellcolor[HTML]{FFFFC7}Relative Unfairness of $F_{dev}$}} &
  $\alpha=5$ &
  0.79 &
  {\color[HTML]{333333}0.48} &
  0.65 &
  0.24 &
  0.73 
\end{tabular}
\\
\vspace{0.2cm}
\begin{tabular}{>{\columncolor[HTML]{FFFFC7}}c l|
>{\columncolor[HTML]{67FD9A}}c c
>{\columncolor[HTML]{67FD9A}}c >{\columncolor[HTML]{ECF4FF}}c }
\multicolumn{2}{c|}{\cellcolor[HTML]{FFFFC7}Model}  &
  {\color[HTML]{333333}$M_{m-s,5}$} &
  $M_{d-s,3}$ &
%  $M_{d-s,3(B)}$ &
  $M_{d-s,5} $&
  $\text{Original}$ 
  \\ \hline
\multicolumn{1}{c|}{\cellcolor[HTML]{FFFFC7}} &
  $\alpha=0$ &
  {\color[HTML]{333333}1.06} &
  0.94 &
 % 0.99 &
  1.00 &
  1 
\\ \cline{2-6} 
\multicolumn{1}{c|}{\cellcolor[HTML]{FFFFC7}} &
  $\alpha=1$ &
  {\color[HTML]{333333}0.86} &
  0.89 &
 % 0.89 &
  0.81 &
  1 \\ \cline{2-6} 
\multicolumn{1}{c|}{\cellcolor[HTML]{FFFFC7}} &
  $\alpha=2$ &
  {\color[HTML]{333333}0.69} &
  0.84 &
 % 0.81 &
  0.66 &
  1 
  \\ \cline{2-6} 
\multicolumn{1}{c|}{\multirow{-4}{*}{\cellcolor[HTML]{FFFFC7}Relative Unfairness of $F_{dev}$}}&
  $\alpha=5$ &
  {\color[HTML]{333333}0.43} &
  0.71 &
 % 0.69 &
  0.44 &
  1 
\end{tabular}
\caption{Relative Unfairness of Deviation-weighted Fairness in the Original and $K$-segment Models.}
\label{table3.5}
\end{table}

% The $K$-segment models with larger $K$ values perform better in terms of deviation-weighted fairness. 
% For $\alpha=2$, the unsmoothed 3-segment model's $RU_{dev}$ value is 0.90, marking a 10\% improvement compared to the original model. The unsmoothed 5-segment model showcases a $RU_{dev}$ value of 0.74, achieving a 26\% improvement over the original model, highlighting the benefits of increasing $K$.

% Smoothing methods substantially enhance the deviation-weighted fairness in $K$-segment models. Consider the case with $\alpha=2$. The Quantile-based method achieves best performance in both cases. When $K=3$, the Quantile-based method lowers the $RU_{dev}$ value to 0.78, which signifies a 13\% reduction  after smoothing. When $K=5$, smoothing reduces the $RU_{dev}$ value to 0.43 -- a significant 42\% reduction.

% As previously mentioned, the value of $\alpha$ affects the sensitivity of $F_{dev}$, with $\alpha=0$ being an exceptional case. Nevertheless, for other values of $\alpha$, the previously noted trends regarding $K$ and smoothing methods hold. Notably, across chosen $\alpha$ values, the Quantile-based 5-segment model consistently performs best in terms of deviation-weighted fairness. The Quantile-based method consistently outperforms others in deviation-weighted fairness across both $K = 3$ and $K = 5$ scenarios.

% This section delivers the following key insights:
\begin{enumerate}
    \item Larger $K$ values enhance deviation-weighted fairness, with 5-segment models generally outperforming 3-segment ones;
    \item The implementation of a suitable smoothing method significantly enhances fairness. The quantile-based method proves to be the most effective for both $K=3$ and $K=5$.
  %  \item Larger $\alpha$ values lead to increased deviation-weighted fairness scores and a reduction in relative unfairness ($RU_{dev}$). While the absolute scores change, the relative performance of the models in comparison to each other remains consistent for non-zero $\alpha$ values.
\end{enumerate}

\subsection{Accuracy--Fairness Analysis across Models}\label{se3}

This section analyzes the trade-off between accuracy and fairness in our models. Figure~\ref{AF} shows a scatter plot of the accuracy ($R^2$) versus the two fairness measures and the Pareto frontier of the various models. These models are differentiated using color coding: blue for the original model, red for the 3-segment model, and green for the 5-segment model.
%\footnote{For clarity in Figure~\ref{AF}, we abbreviate the model $M_{\text{smoothing method},K}$ notation to simply `\text{smoothing method},K'.} 
The blue line represents the accuracy--fairness frontier, forming a convex hull that connects models excelling along both the accuracy and fairness dimensions. Here are a few observations in order:

\begin{figure}[H]
\centering  
\subfigure[Accuracy ($R^2$) versus $F_{grp}$]{
\label{af1}
\includegraphics[width=0.45\textwidth]{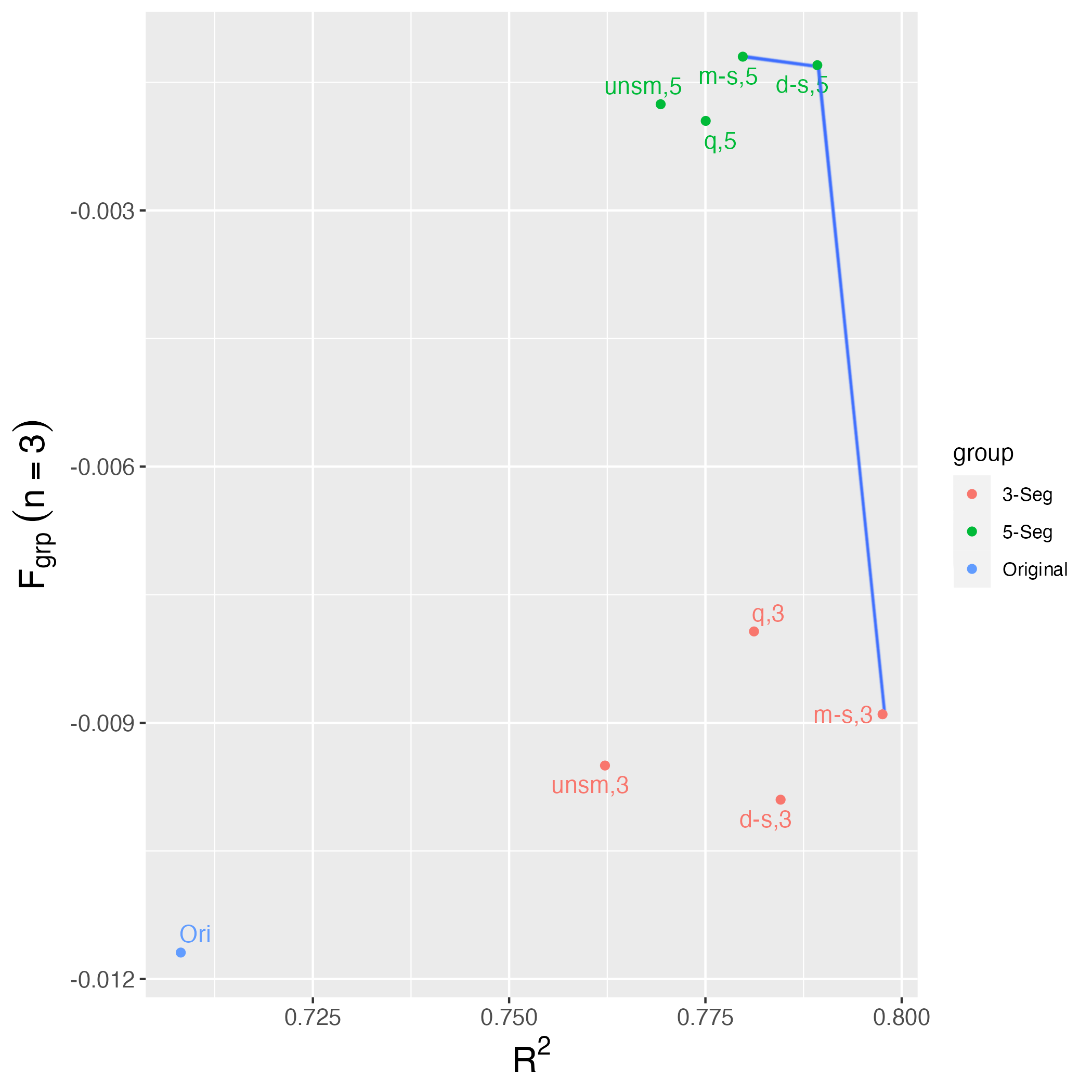}}
\subfigure[Accuracy ($R^2$) versus $F_{dev}$]{
\label{af2}
\includegraphics[width=0.45\textwidth]{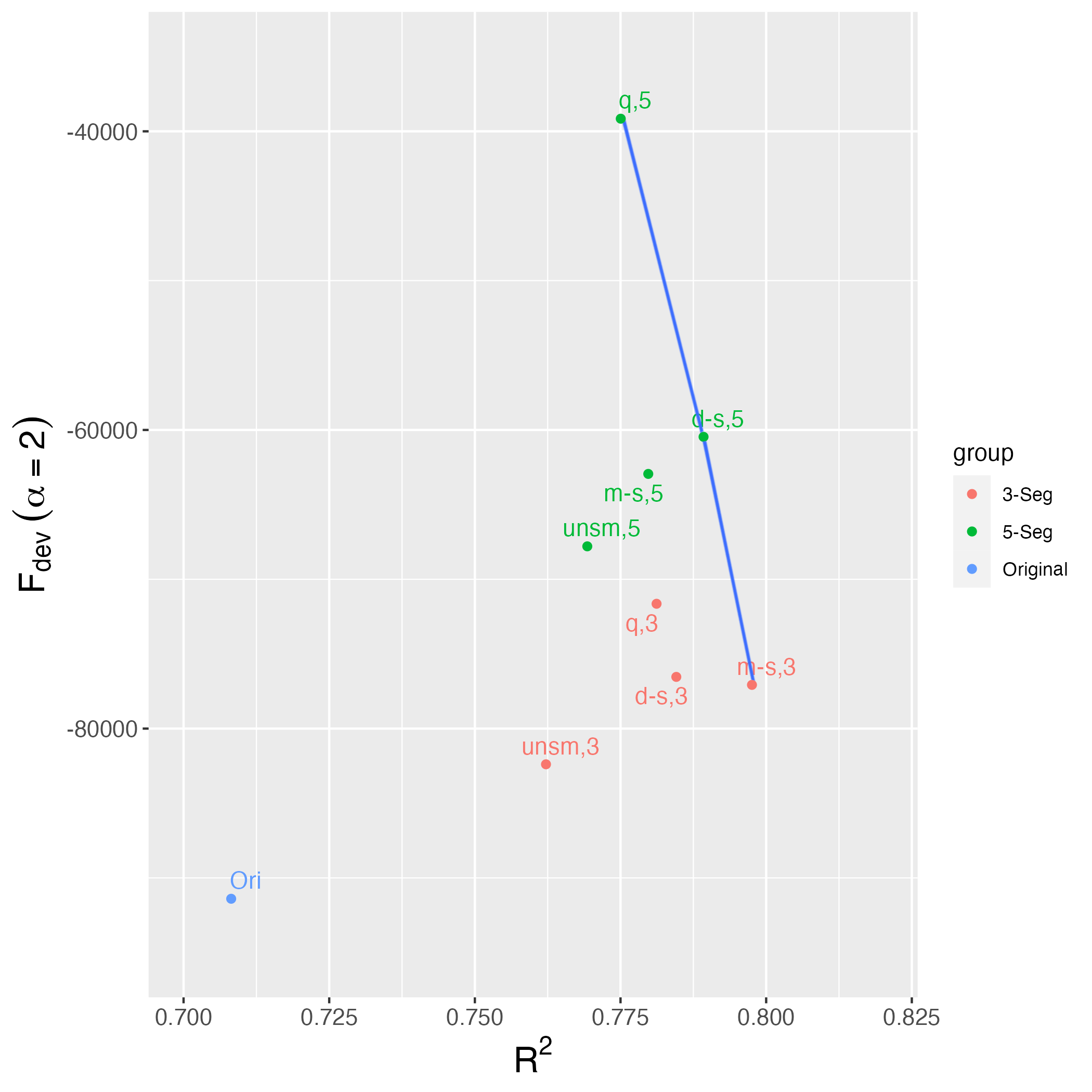}}
\caption{Analysis of the Accuracy--Fairness Trade-off across Fairness Measures.}
\label{AF}
\end{figure}

\begin{enumerate}
    \item The $K$-segment model improves both the accuracy ($R^2$) and fairness measures compared to the original model. In general, the model is fairer when $K=5$ than when $K=3$, while keeping a similar level of accuracy. Smoothing methods can improve both the accuracy and fairness measures.
    \item Models $M_{\text{d-s},5}$ and $M_{\text{m-s},3}$ appear in the frontier, which showcases their superior performance in terms of both accuracy and fairness.
    \item There is a clear cluster pattern in Figure \ref{AF}(a), but not in Figure \ref{AF}(b). This is due to the fundamental difference between the two fairness measures.
    
    % A distinct `cluster' pattern emerges, where models with the same $K$ value group closely together in the accuracy--fairness space. Segmentation enhances both accuracy and fairness measures, particularly for group fairness. The improvement is noticeable at $K=3$ and becomes even more significant at $K=5$.
    % \item The group fairness and deviation-weighted fairness measures, while related, shows different aspects of model and data. 
    % \item Smoothing methods show potential to boost both accuracy and fairness, though their effectiveness varies based on the specific method and parameters used.
    % \item The selected model for the Assessment set, $M_{\text{d-s},5}$, balances accuracy and fairness. It significantly outperforms another frontier model, $M_{\text{m-s},3}$, in fairness, despite a minor trade-off in accuracy. 
\end{enumerate}

\section{Conclusions and Future Directions}\label{con}

The pervasive regressivity evident in our empirical data motivated us to develop the $K$-segment model.
This model works by first splitting the data into multiple disjoint subsets based on the previous year's property value predictions and then conducting estimates within each group separately.
Intuitively, such a split can be obtained by focusing on different quantiles of the previous year's predictions. However, doing so creates discontinuities in the assessed values that make the model harder to interpret in practice and harder to justify. To address this issue, 
we also introduce a set of different smoothing methods (quantile-based smoothing, midpoint score-based smoothing, and distance score-based smoothing). 
 Our results indicate that a 5-segment model, in conjunction with the distance score-based smoothing method, substantially improves both accuracy and fairness (in particular group fairness and deviation-weighted fairness) relative to the model employed in practice.

Our findings suggest several avenues for future research. (1) There is a need to analyze the temporal dynamics of the $K$-segment model by evaluating its performance over different years. The performance and fairness measures of the current model are based primarily on data for 2022, but factors such as interest rates, economic growth, population trends, and policy changes can result in substantial year-to-year variations. These variations can affect the distribution of data and, consequently, the model's performance.
(2) The potential impact of gentrification, regional and demographic development, and tax policy on our fairness measures should not be overlooked. These factors could create new tax inequalities, which might not be captured by the fairness measures derived from a single year's data. Thus, investigating the robustness of the $K$-segment model under these dynamic conditions can provide insights into the model's adaptability and generalizability.
% (3) Extending the application of the $K$-segment model beyond the geographical region we focus on in this study could yield invaluable insight into the model’s performance with property datasets from different regions, or even with other types of data that are prone to regressivity.
(3) The versatility of the model in addressing regressivity and promoting fairness should be explored in other domains such as income taxation, healthcare, and education. This expansion of applications can demonstrate the broader applicability of the $K$-segment model.
(4) Even though our model's generalizability has been bolstered through cross-validation techniques in each submodel, other robustness checks, such as bootstrap sampling or different machine learning techniques, could be valuable. Assuring robustness in the face of diverse statistical and practical challenges is crucial for the effective real-world application of the $K$-segment model.

\bibliographystyle{IEEEtran}
\bibliography{sample}

\newpage
\appendix
\renewcommand{\thesection}{\Alph{section}}

\section{Supplementary Plots and Visual Data Analysis}\label{app4}

\subsection{Plots for the Unsmoothed $K$-segment Model}\label{app1}

\begin{figure}[H]
\centering  
\subfigure[$r$ versus $\log(x)$, $K=3$]{
\includegraphics[width=0.3\textwidth]{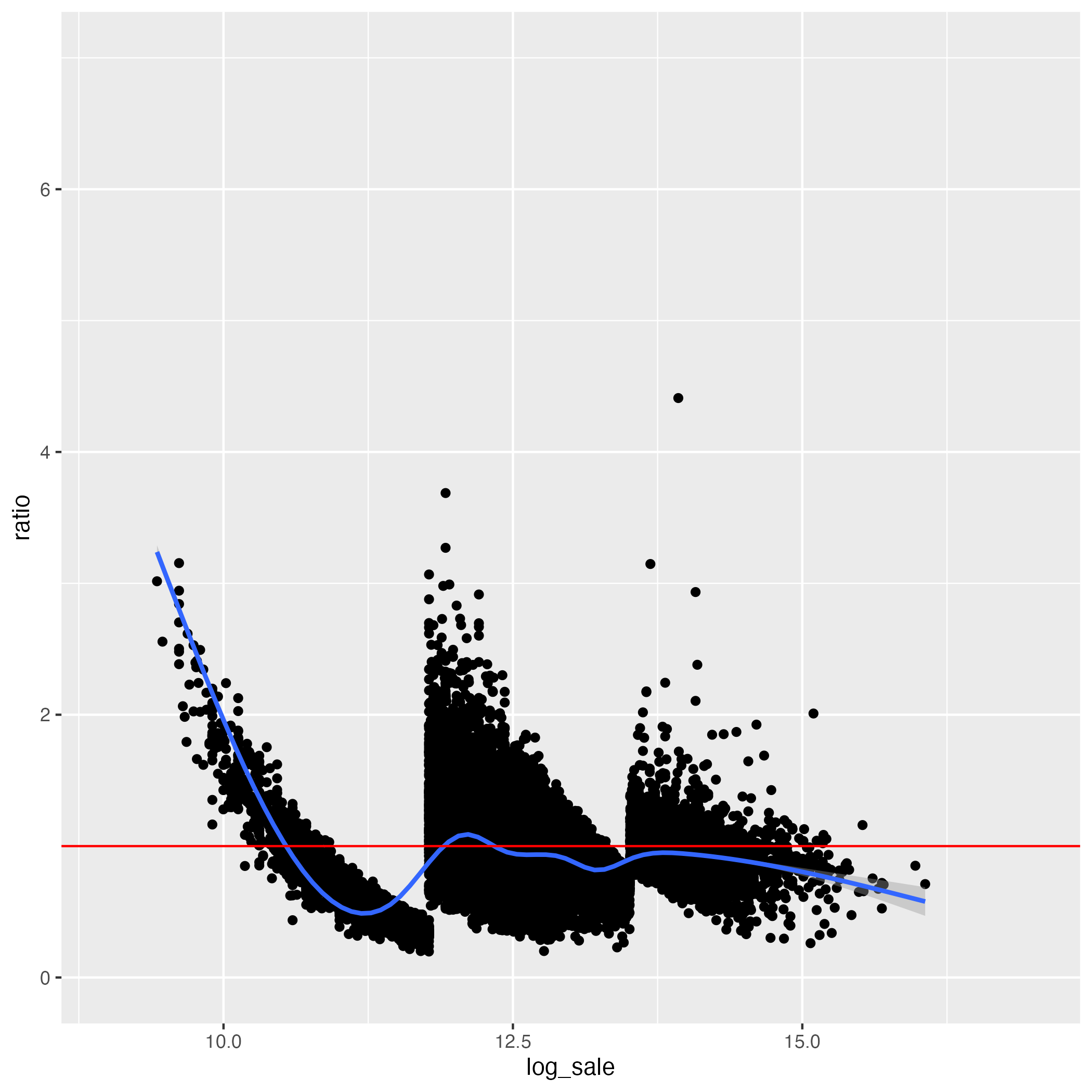}}
% \subfigure[$\log(V)$ versus $x$, $K=3$]{
% \includegraphics[width=0.3\textwidth]{plot10/unsmooth/unsm3ori_logratio_test.png}}
% \subfigure[Trend Scatter Plot of $Q$, $K=3$]{
% \includegraphics[width=0.3\textwidth]{plot10/unsmooth/unsm3s_div_ass.png}}
\subfigure[$r$ versus $\log(x)$, $K=5$]{
\includegraphics[width=0.3\textwidth]{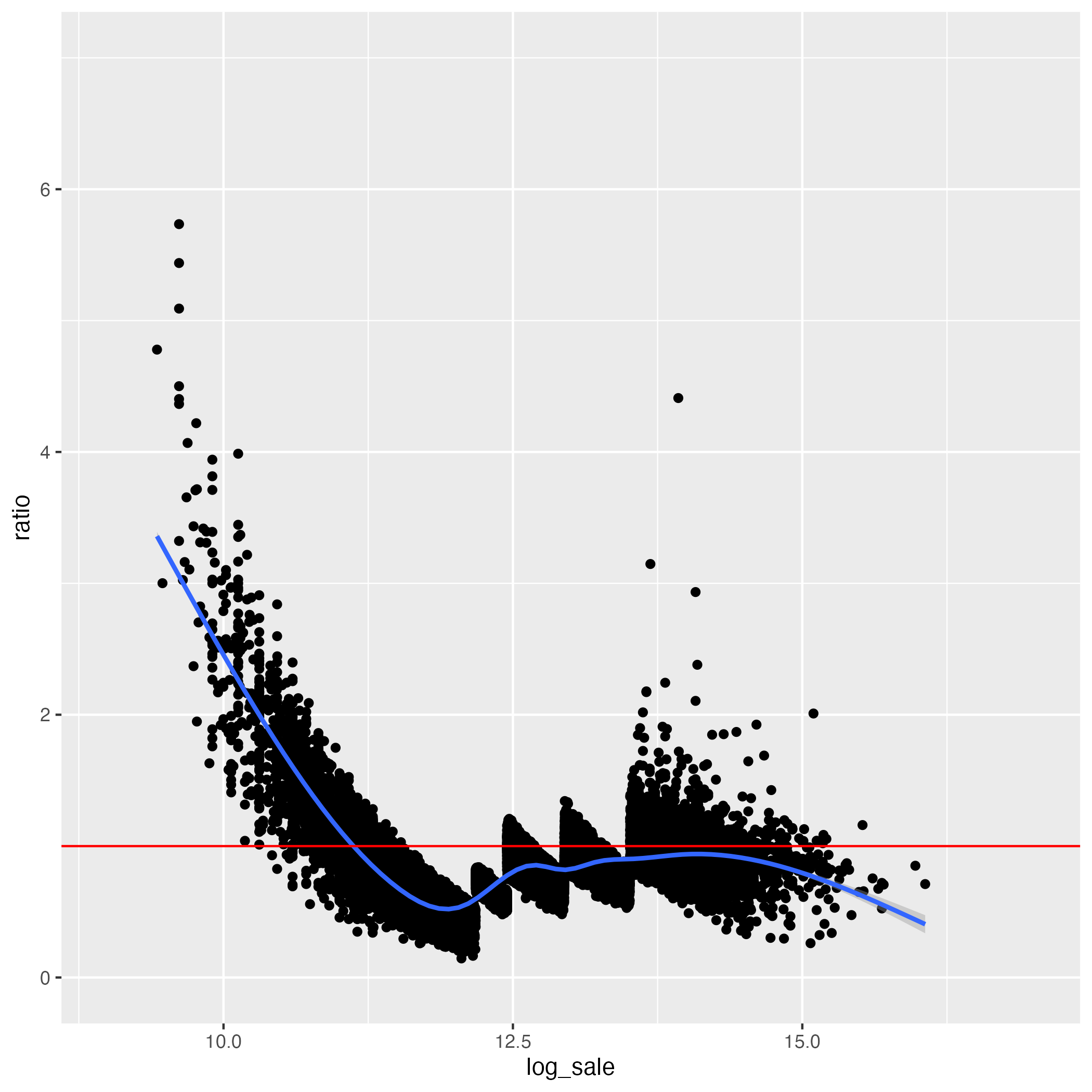}}
% \subfigure[$\log(V)$ versus $x$, $K=5$]{
% \includegraphics[width=0.3\textwidth]{plot10/unsmooth/unsm5ori_log_test.png}}
% \subfigure[Trend Scatter Plot of $Q$, $K=5$]{
% \includegraphics[width=0.3\textwidth]{plot10/unsmooth/unsm5s_div_ass.png}}
\subfigure[Comparative Bar Plot, $K=3$]{
\includegraphics[width=0.45\textwidth]{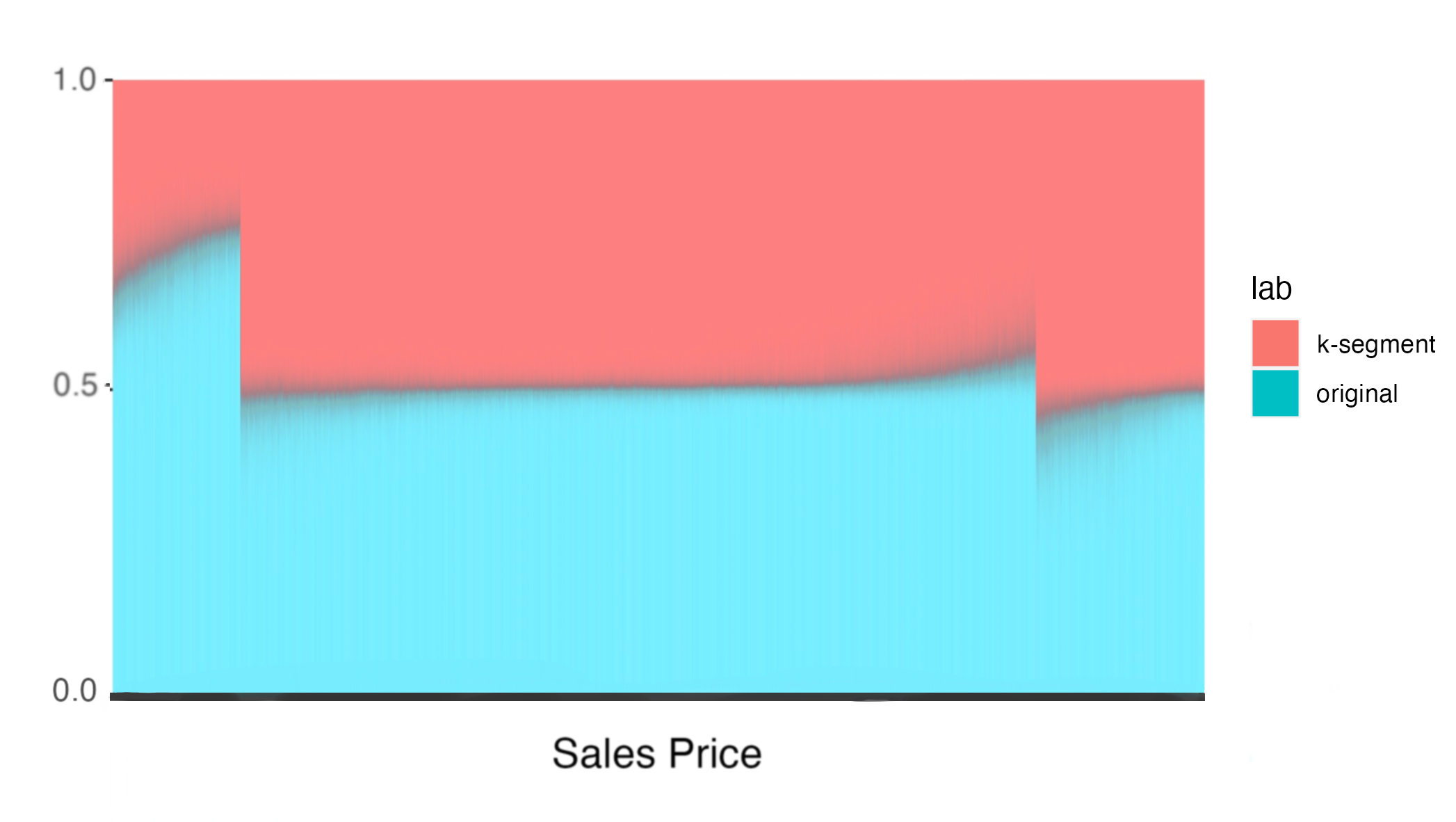}}
\subfigure[Comparative Bar Plot, $K=5$]{
\includegraphics[width=0.45\textwidth]{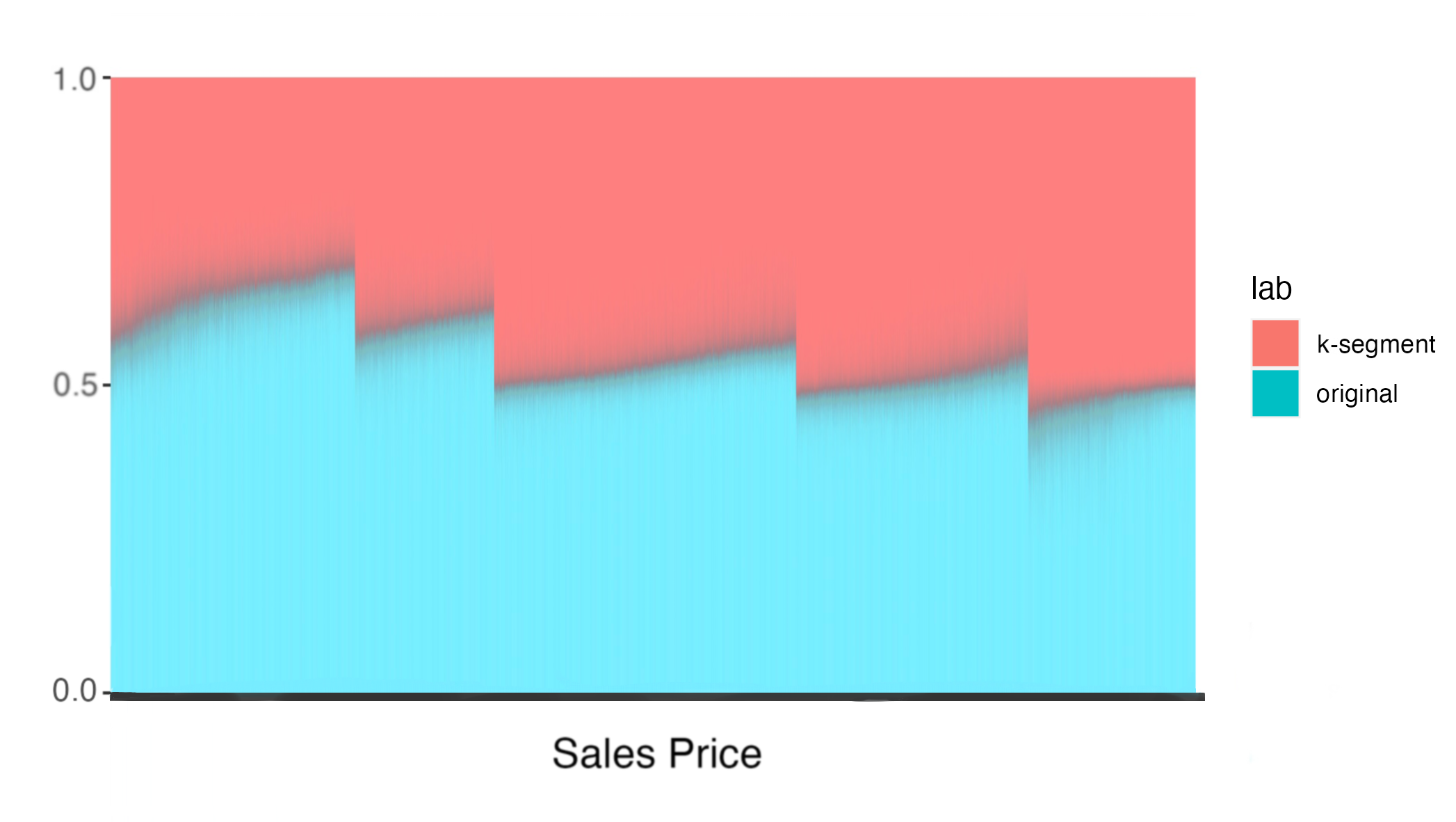}}
\caption{Supplementary Analysis: Behavior of the Unsmoothed $K$-segment Model.}
\label{fig:app1}
\end{figure}

Figure~\ref{fig:app1} provides supplementary illustrations of the unsmoothed $K$-segment model for $K=3$ and $K=5$, which emphasize the conspicuous and undesirable jumps between segments in the unsmoothed $K$-segment model.

Figure~\ref{fig:app1} shows the effectiveness of the unsmoothed $K$-segment model in addressing regressivity, and compares it with the $K$-segment models that incorporate smoothing methods, as discussed in Section~\ref{modelfor}.
Specifically, panels (a) and (b) in Figure~\ref{fig:app1} depict the relationship between ratio $r$ and logarithm of sales price $\log x$, refer to Figure~\ref{ratio2} for a similar illustration. Panels (c) and (d) show the comparative illustration of the unsmoothed $K$-segment and original models by sales price, see Figure~\ref{per} for a similar illustration.

\subsection{High Precision of the $K$-segment Model}

\begin{figure}[H]
\centering  
\subfigure[Original Model]{
\includegraphics[width=0.18\textwidth]{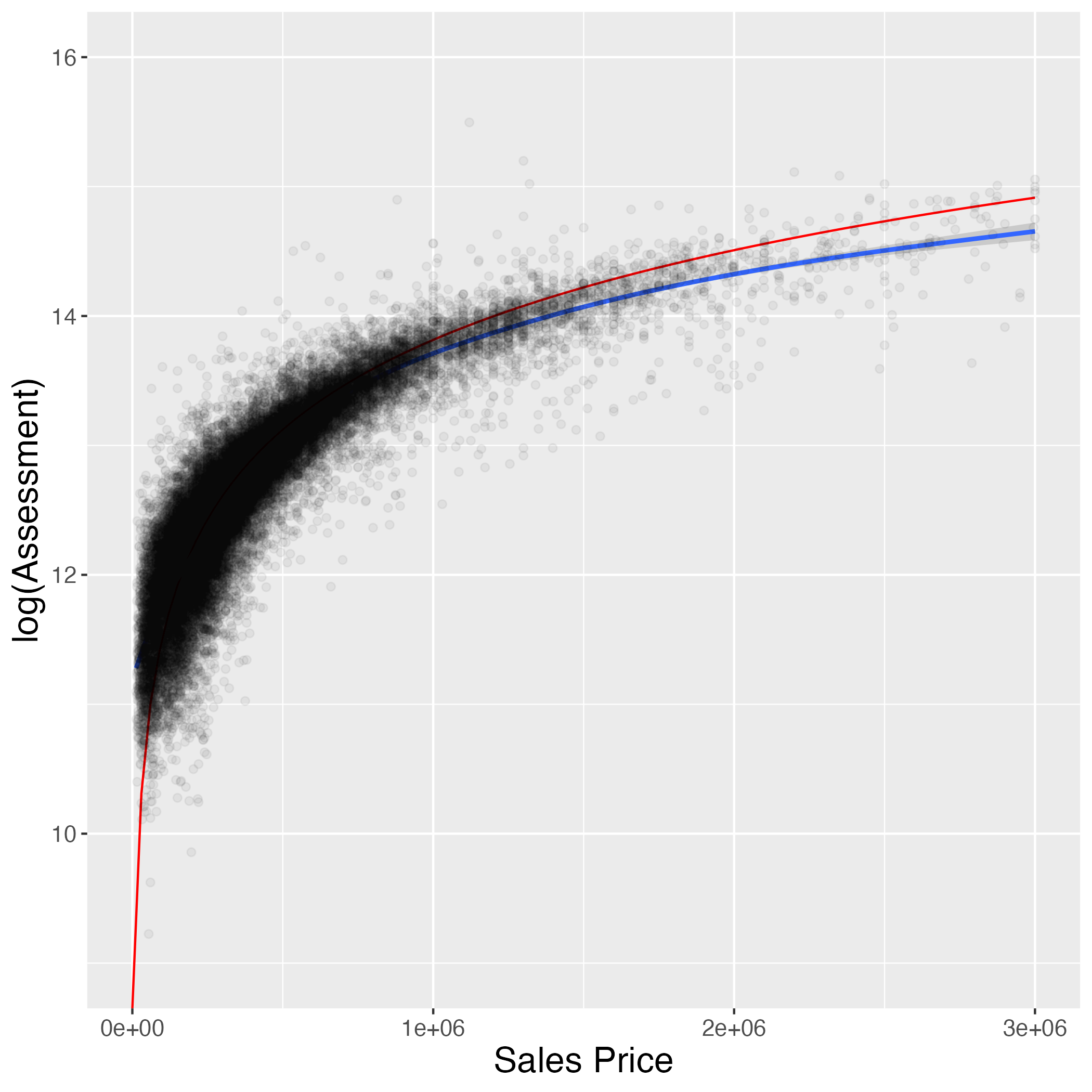}}
\subfigure[$\textit{M}_{\textit{d--s},\textit{3}}$]{
\includegraphics[width=0.18\textwidth]{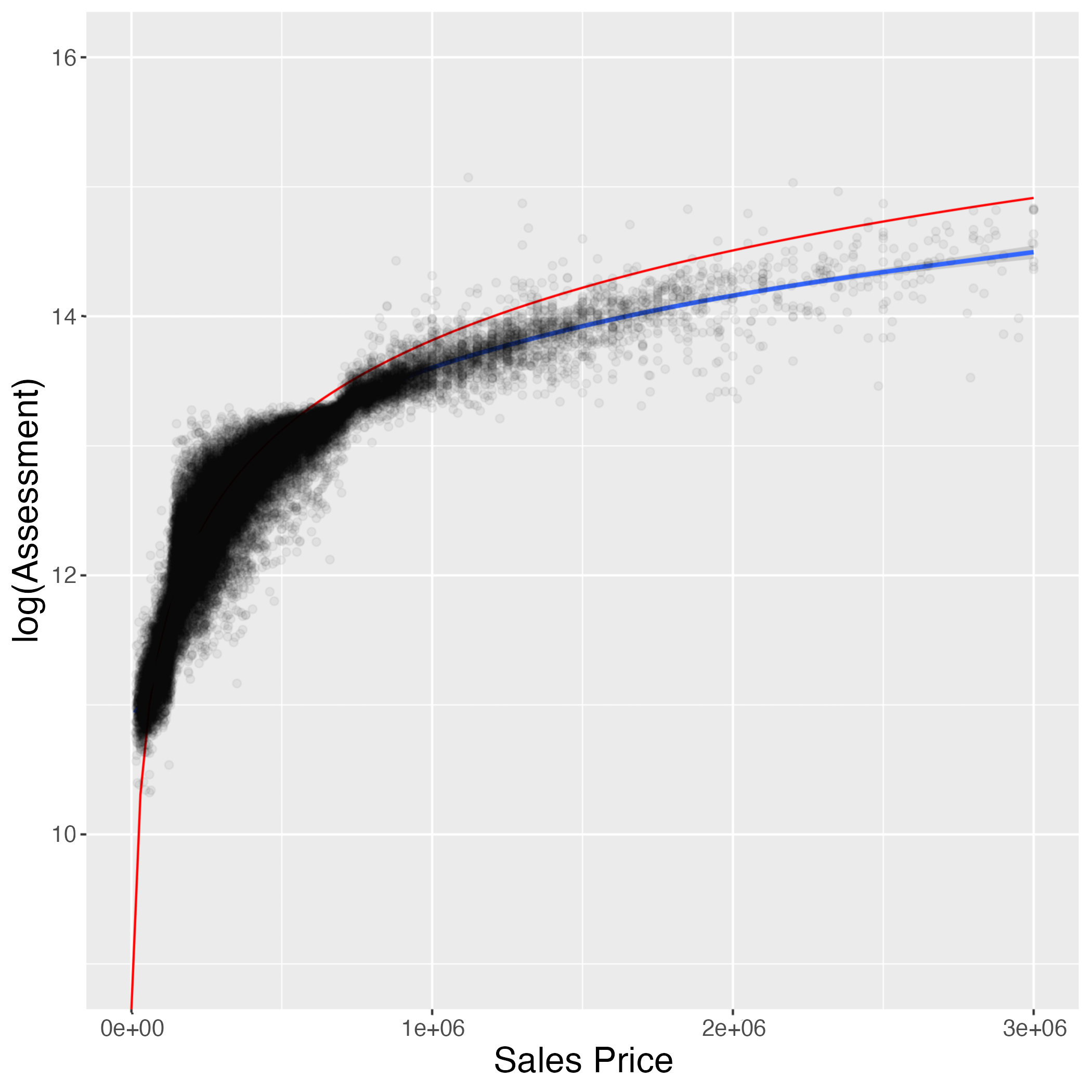}}
\subfigure[$\textit{M}_{\textit{q},\textit{5}}$]{
\includegraphics[width=0.18\textwidth]{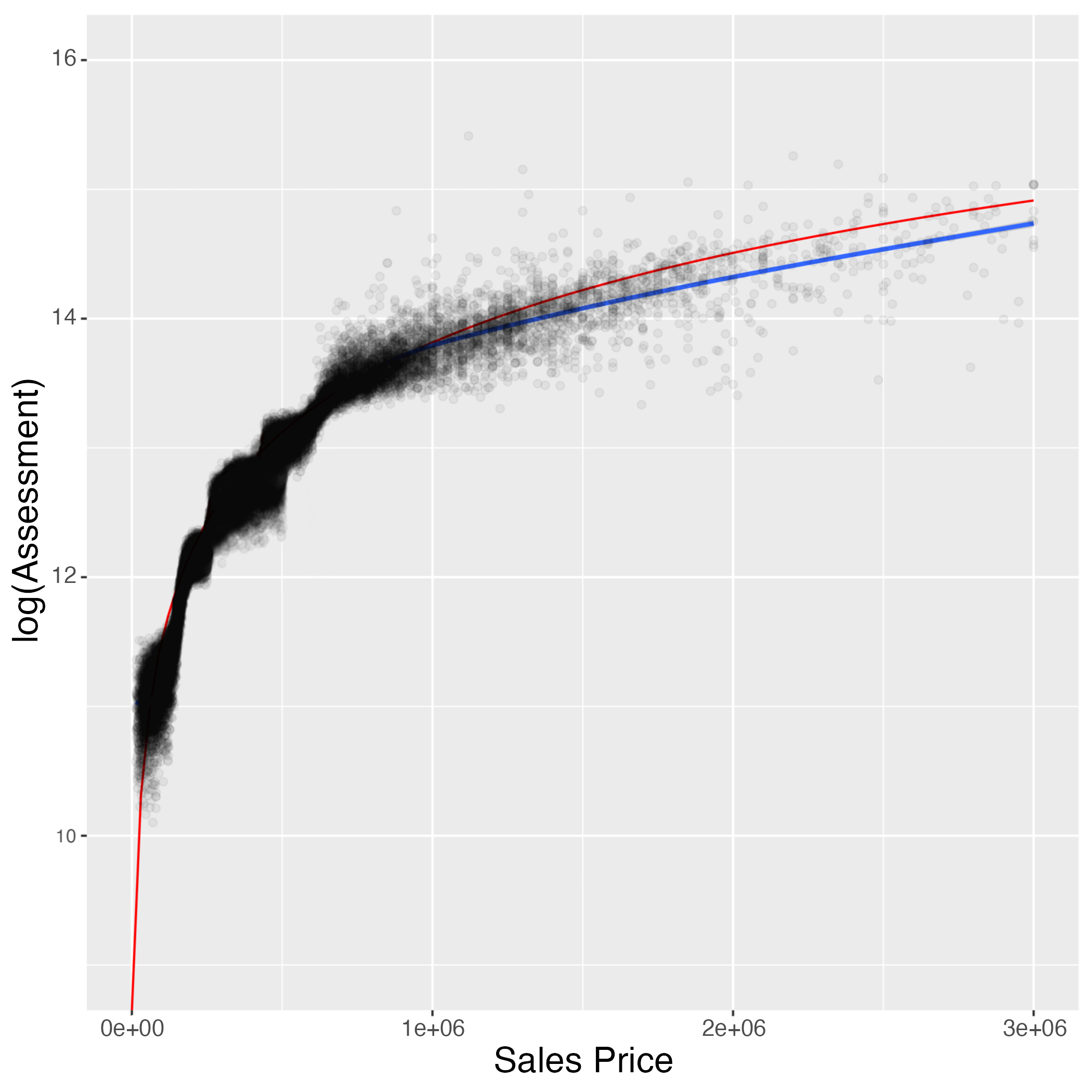}}
\subfigure[$\textit{M}_{\textit{m--s},\textit{5}}$]{
\includegraphics[width=0.18\textwidth]{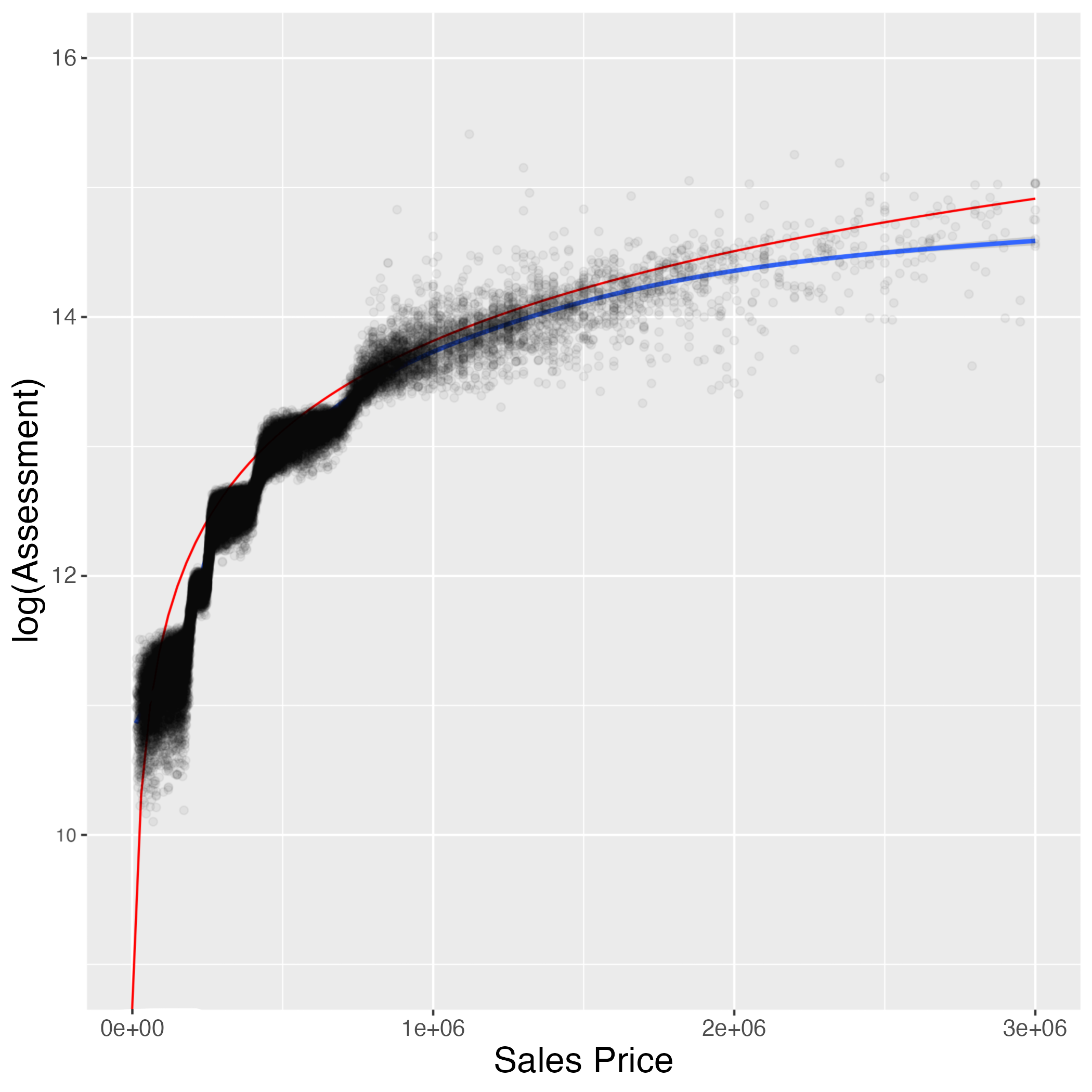}}
\subfigure[$\textit{M}_{\textit{d--s},\textit{5}}$]{
\includegraphics[width=0.18\textwidth]{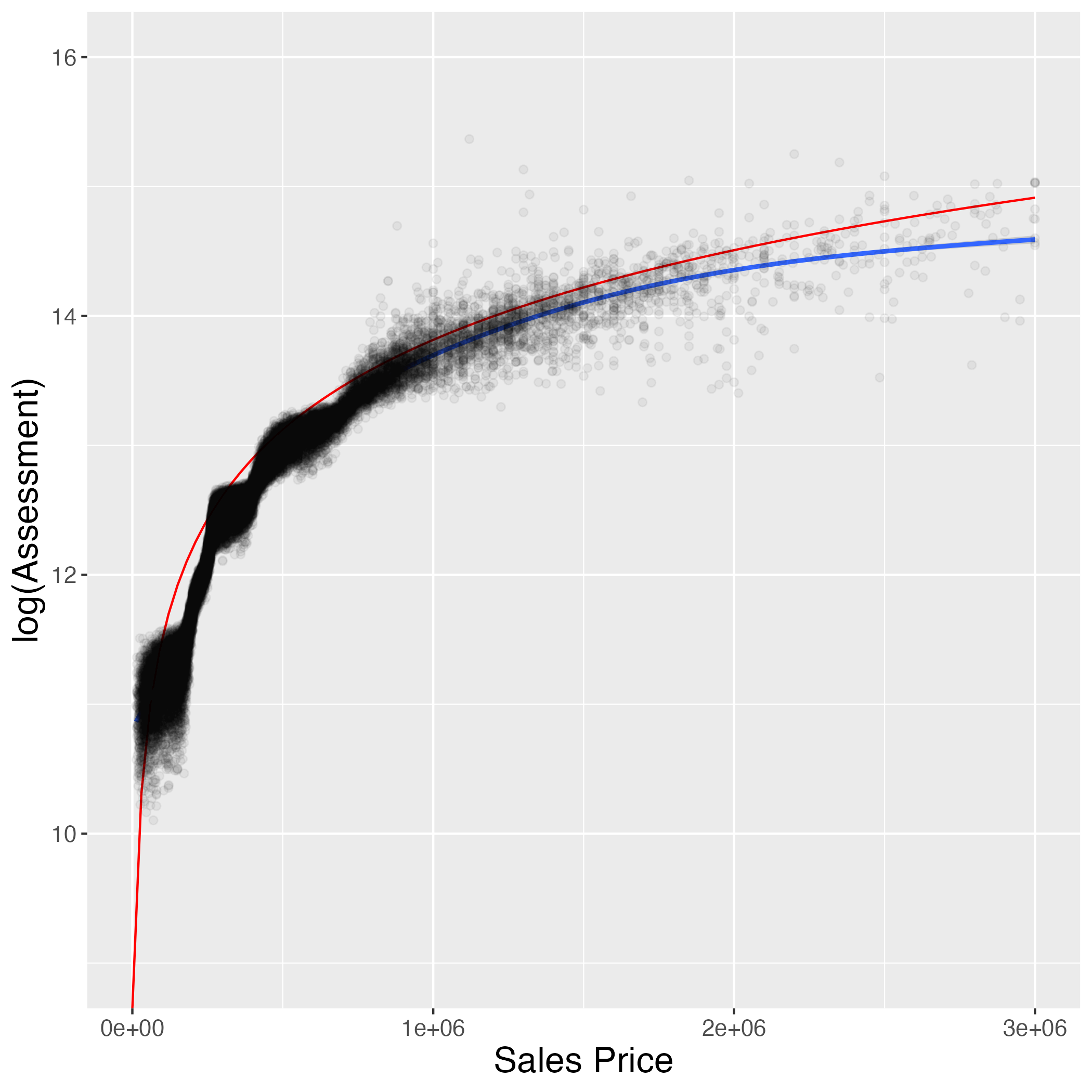}}
\caption{Logarithm of Assessment ($\log (v)$) versus Sales Price ($x$).}\label{11b}
\end{figure}

In Figure~\ref{11b}, the x-axis is the sales price and the y-axis is the logarithm of the assessment value. The red line represents the baseline, where the assessment is exactly equal to the sales price. 
Figure~\ref{11b} shows that the $K$ segment model shows a higher valuation precision compared to the original model for properties sold in 2021, and thus smaller random errors.

\subsection{Plots for Selection of Parameter $\mu$}\label{app2}

\begin{figure}[H]
\centering  
\subfigure[$\mu=100$]{
\includegraphics[width=0.45\textwidth]{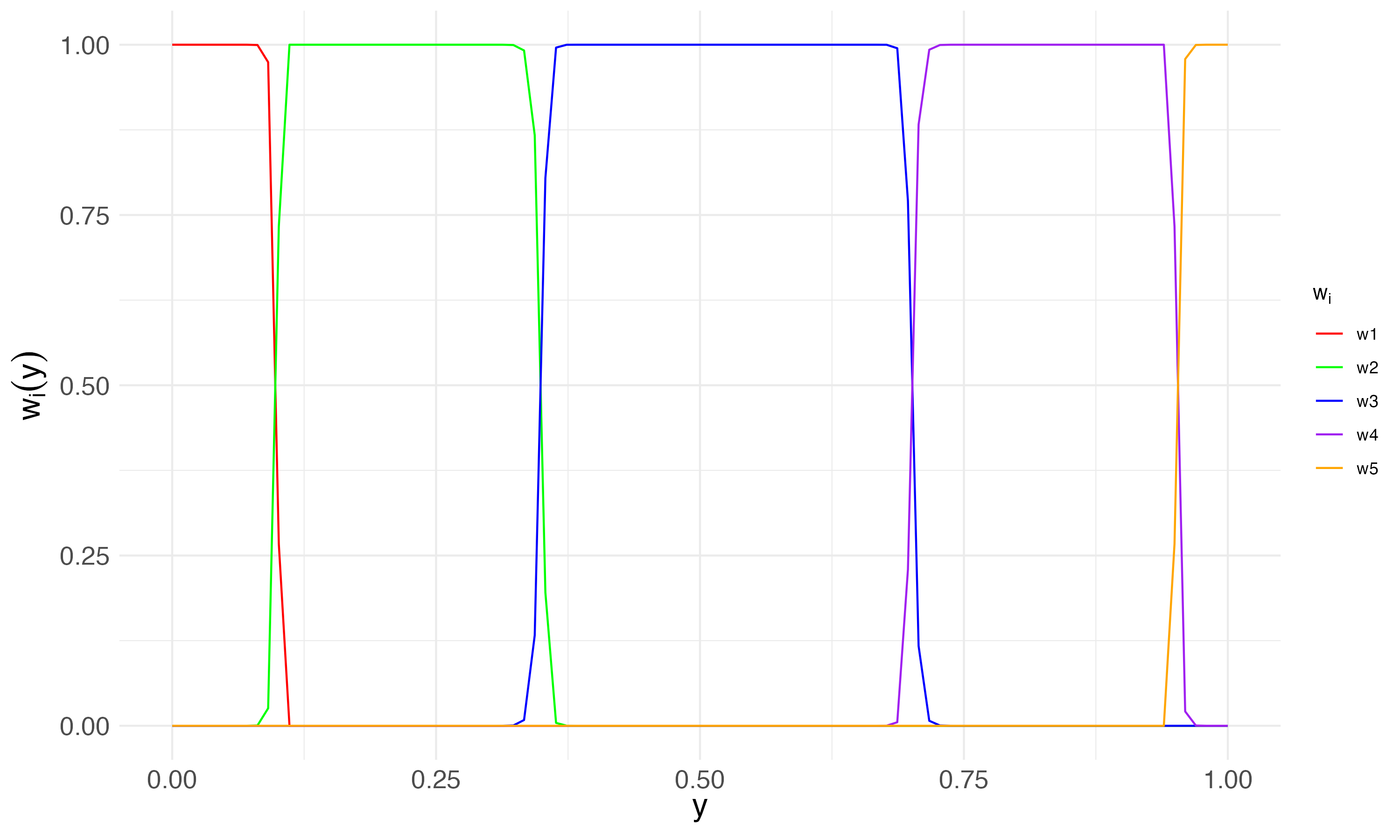}}
\subfigure[$\mu=1$]{
\includegraphics[width=0.45\textwidth]{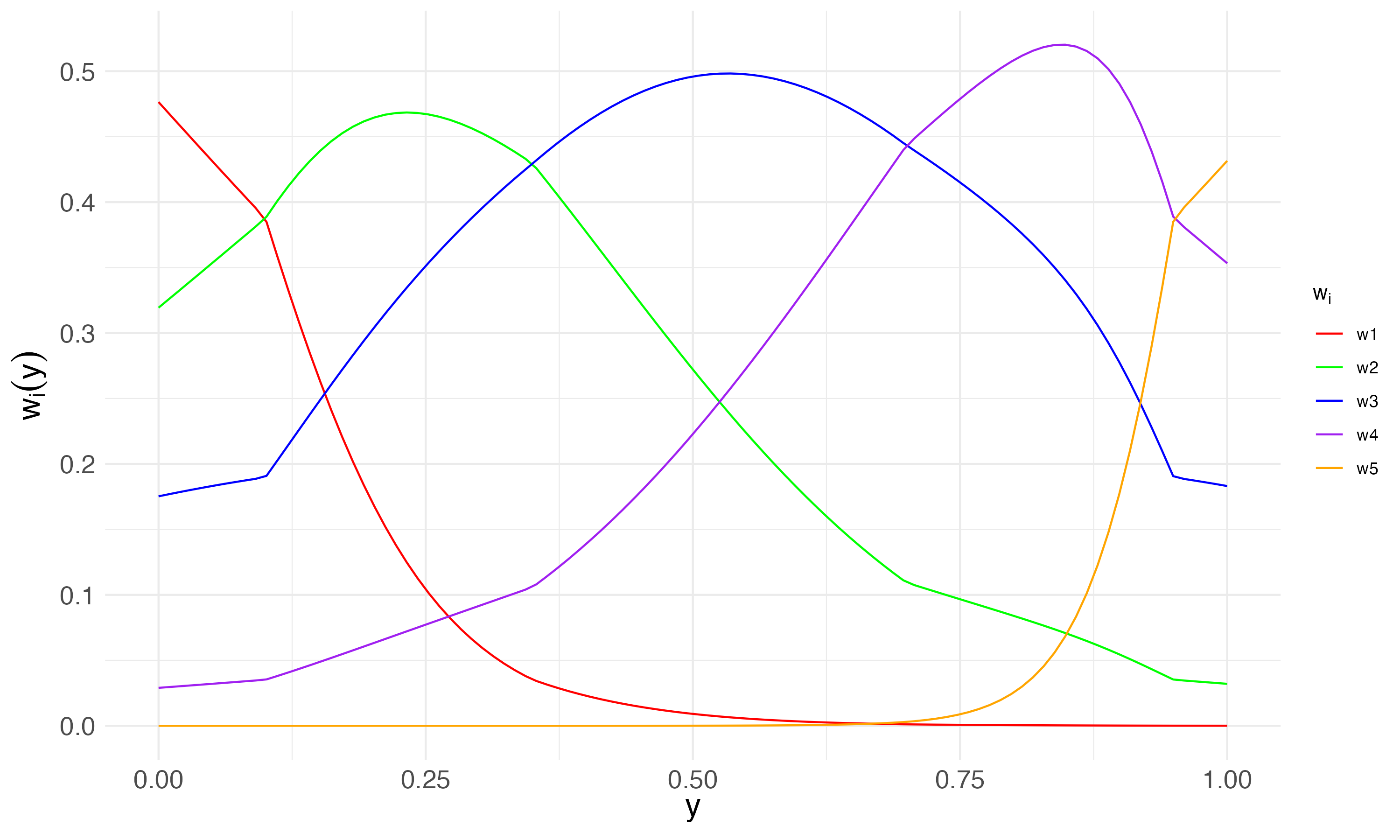}}
\caption{Weighting function $w_i(y)$ for $\textit{M}_{\textit{d--s},\textit{5}}$ versus quantile $y$ (when 
$K=5,\eta_1=0.1,\eta_2=0.35,\eta_3=0.7,\eta_4=0.95$).}
%$K=5,\eta_1=0.2,\eta_2=0.35,\eta_3=0.7,\eta_4=0.9$).}
\label{fig:app4}
\end{figure}

To illustrate the effects of varying the $\mu$ parameter in the score-based smoothing method of the $K$-segment model, we provide appendix plots. Figure~\ref{fig:app4} presents the relationship between the quantile \( y \) and the respective weights \( w_k \) for \( k=1,\ldots,5 \), with panels (a) and (b) corresponding to \(\mu=100\) and \(\mu=1\), respectively.

% \begin{figure}[H]
% \centering  
% \subfigure[$w_1(y)$ against  $y$]{
% \includegraphics[width=0.3\textwidth]{plot10/selectmu/30p1.png}}
% \subfigure[$w_2(y)$ against  $y$]{
% \includegraphics[width=0.3\textwidth]{plot10/selectmu/30p2.png}}
% \subfigure[$w_3(y)$ against  $y$]{
% \includegraphics[width=0.3\textwidth]{plot10/selectmu/30p3.png}}
% \subfigure[$w_4(y)$ against  $y$]{
% \includegraphics[width=0.3\textwidth]{plot10/selectmu/30p4.png}}
% \subfigure[$w_5(y)$ against  $y$]{
% \includegraphics[width=0.3\textwidth]{plot10/selectmu/30p5.png}}
% \caption{
% Visualizing the Weight-Function in a 5-segment Model ($\eta_1=0.2,\eta_2=0.35,\eta_3=0.7,\eta_4=0.9,\mu=30$).}
% \label{fig:app3}
% % \end{figure}

% \begin{figure}[H]
% \centering  
% \subfigure[$w_1(y)$ against  $y$]{
% \includegraphics[width=0.3\textwidth]{plot10/selectmu/1p1.png}}
% \subfigure[$w_2(y)$ against  $y$]{
% \includegraphics[width=0.3\textwidth]{plot10/selectmu/1p2.png}}
% \subfigure[$w_3(y)$ against  $y$]{
% \includegraphics[width=0.3\textwidth]{plot10/selectmu/1p3.png}}
% \subfigure[$w_4(y)$ against  $y$]{
% \includegraphics[width=0.3\textwidth]{plot10/selectmu/1p4.png}}
% \subfigure[$w_5(y)$ against  $y$]{
% \includegraphics[width=0.3\textwidth]{plot10/selectmu/1p5.png}}

Figure~\ref{k5} provides the $\mu=10$ case. When $\mu$ is too large, the efficacy of the weighting function diminishes; conversely, a too small $\mu$ fails to adequately represent the emphasis on a specific segment.

% \subsection{Sales Price Distribution}
% \begin{figure}[H]
%     \centering
%     \includegraphics[width=0.6\textwidth]{plot10/compare/distri_sale price.png}
%     \caption{Illustration of Data Skewness in Property Assessments.}
%     \label{skew}
% \end{figure}

% In Figure~\ref{skew}, the y-axis represents the count of properties at each price, and the x-axis represents the sales price.
% We can see that the distribution of sales prices is indeed skewed to the right. This is a common characteristic of real estate markets, where there are many more lower-priced homes.
% This skewness in the data could potentially affect the performance of the property assessment models. Models trained on such data might learn to predict more accurately for the majority class; in this case, it might exacerbate regressivity by overvaluing lower-priced homes and undervaluing higher-priced ones.

% \begin{figure}[H]
%     \centering
%     \includegraphics[width=0.6\textwidth]{plot10/compare/compare 2 model.png}
%     \caption{Comparison of the StA ratios for the Original Model and the $M_{\text{d-s},5}$ Model.}
%     \label{compare}
% \end{figure}

% In Figure~\ref{compare}, we present a comparative visualization of the StA ratio for the original model and the $K$-segment model, $M_{\text{d-s},5}$. The x-axis denotes the sales price, while the y-axis corresponds to the StA ratio. The red dashed line shows an ideal StA ratio of 1, indicating a scenario where the valuation is perfectly equivalent to the sales price.

\section{Actual Values of Fairness Measures}

This appendix presents the raw values of $F_{grp}$ and $F_{dev}$, independent of any associated unfairness definitions.

\subsection{Group Fairness}\label{app3}
\begin{table}[H]
\centering
\setlength{\tabcolsep}{11mm}{
\begin{tabular}{c|cc}
\rowcolor[HTML]{FFFFC7} 
\cellcolor[HTML]{FFFFC7}                        & \multicolumn{2}{c}{\cellcolor[HTML]{FFFFC7}Group Fairness}             \\ \cline{2-3} 
\multirow{-2}{*}{\cellcolor[HTML]{FFFFC7}Model} & \multicolumn{1}{c|}{$n=2$}                                 & $n=3$         \\ \hline
$M_{unsm,3}$                            & \multicolumn{1}{c|}{-0.00217628}                          & -0.009524356 \\
\rowcolor[HTML]{67FD9A} 
$M_{unsm,5}$         & \multicolumn{1}{c|}{\cellcolor[HTML]{67FD9A}-0.0004896838} & -0.001756903 \\
$M_{q,3}$                    & \multicolumn{1}{c|}{-0.001840958}                         & -0.007934671 \\
\rowcolor[HTML]{67FD9A} 
$M_{q,5}$ & \multicolumn{1}{c|}{\cellcolor[HTML]{67FD9A}-0.0003966192} & -0.001952189 \\
$M_{m-s,3}$                      & \multicolumn{1}{c|}{-0.002049917}                         & -0.008906032 \\
\rowcolor[HTML]{67FD9A} 
$M_{m-s,5}$   & \multicolumn{1}{c|}{\cellcolor[HTML]{67FD9A}-0.0003587153} & -0.001217768 \\
$M_{d-s,3}$                           & \multicolumn{1}{c|}{-0.002257118}                         & -0.00988218  \\
%$M_{d-s,3(B)}$                       & \multicolumn{1}{c|}{-0.001920633}                         & -0.008014704 \\
\rowcolor[HTML]{67FD9A} 
$M_{d-s,5}$        & \multicolumn{1}{c|}{\cellcolor[HTML]{67FD9A}-0.0003408842} & -0.001286112 \\
\rowcolor[HTML]{ECF4FF} 
Original                              & \multicolumn{1}{c|}{\cellcolor[HTML]{ECF4FF}-0.002612381} & -0.01168704 
\end{tabular}}
\caption{Group Fairness ($F_{grp}$) in the Original and $K$-segment Models.}
\label{tab:table2}
\end{table}

% \begin{remark}
%     [Approaches to tackling large-scale challenges in managing data of 550k by 550k]
%      Pairing each data point with another for processing is a formidable challenge. Particularly, when $n$ equals 2, this would imply pairing 550k data points with another 550k data points. A for-loop implementation for this task would be impractical, as computing a single pair of differences every second would take a year to complete. Even using matrices would be infeasible, as the memory required to handle a 550k$\times$550k matrix would be excessive. To tackle this issue, we adopt the chunking method, which breaks the computation into smaller segments.
% \end{remark}

As can be seen in Table~\ref{tab:table2}, 
the $K$-segment models enhance the fairness by increasing the $F_{grp}$ values compared to those in the original model.
However, the absolute magnitudes of the $F_{grp}$ values are notably low. To create a more nuanced and meaningful analysis, we use the $RU_{grp}$ measure, detailed in Section~\ref{sectiongrp}, in subsequent discussions.
For our experiments, we set $\mu = 10$ in the score-based smoothing method.

\subsection{Deviation-weighted Fairness}\label{app3;2}

\begin{table}[H]
\centering
\setlength{\tabcolsep}{4.8mm}{
\begin{tabular}{c|cccc}
\rowcolor[HTML]{FFFFC7} 
\cellcolor[HTML]{FFFFC7}                        & \multicolumn{4}{c}{\cellcolor[HTML]{FFFFC7}Deviation-weighted Fairness}                                                 \\ \cline{2-5} 
\multirow{-2}{*}{\cellcolor[HTML]{FFFFC7}Model} & $\alpha=0$                        & $\alpha=1$                        & $\alpha=2$                        & $\alpha=5$   \\ \hline
$M_{unsm,3}$                          & \multicolumn{1}{c|}{-225094.6} & \multicolumn{1}{c|}{-129406.4} & \multicolumn{1}{c|}{-82391.24} & -33411.61 \\
\rowcolor[HTML]{67FD9A} 
$M_{unsm,5}$ &
  \multicolumn{1}{c|}{\cellcolor[HTML]{67FD9A}-236884.3} &
  \multicolumn{1}{c|}{\cellcolor[HTML]{67FD9A}-121708.6} &
  \multicolumn{1}{c|}{\cellcolor[HTML]{67FD9A}-67794.16} &
  -20176.55 \\
$M_{q,3}$                   & \multicolumn{1}{c|}{-201342.2} & \multicolumn{1}{c|}{-114525.3} & \multicolumn{1}{c|}{-71640.49} & -27270.65 \\
\rowcolor[HTML]{67FD9A} 
$M_{q,5}$ &
  \multicolumn{1}{c|}{\cellcolor[HTML]{67FD9A}-150699.9} &
  \multicolumn{1}{c|}{\cellcolor[HTML]{67FD9A}-73883.58} &
  \multicolumn{1}{c|}{\cellcolor[HTML]{67FD9A}-39154.31} &
  -10028.59 \\
$M_{m-s,3}$                      & \multicolumn{1}{c|}{-209657.8} & \multicolumn{1}{c|}{-121042.1} & \multicolumn{1}{c|}{-77072.57} & -30893.97 \\
\rowcolor[HTML]{67FD9A} 
$M_{m-s,5}$ &
  \multicolumn{1}{c|}{\cellcolor[HTML]{67FD9A}-224581.9} &
  \multicolumn{1}{c|}{\cellcolor[HTML]{67FD9A}-114261.2} &
  \multicolumn{1}{c|}{\cellcolor[HTML]{67FD9A}-62943.33} &
  -18277.92 \\
$M_{d-s,3}$                          & \multicolumn{1}{c|}{-199236.6} & \multicolumn{1}{c|}{-119241.8} & \multicolumn{1}{c|}{-76537.73} & -29879.35 \\
%$M_{d-s,3(B)}$                     & \multicolumn{1}{c|}{-209575.7} & \multicolumn{1}{c|}{-118724.9} & \multicolumn{1}{c|}{-74264.9}  & -28879.56 \\
\rowcolor[HTML]{67FD9A} 
$M_{d-s,5}$ &
  \multicolumn{1}{c|}{\cellcolor[HTML]{67FD9A}-211011.2} &
  \multicolumn{1}{c|}{\cellcolor[HTML]{67FD9A}-108302.6} &
  \multicolumn{1}{c|}{\cellcolor[HTML]{67FD9A}-60460.62} &
  -18458.86 \\
\rowcolor[HTML]{ECF4FF} 
Original &
  \multicolumn{1}{c|}{\cellcolor[HTML]{ECF4FF}-211218.6} &
  \multicolumn{1}{c|}{\cellcolor[HTML]{ECF4FF}-133576.9} &
  \multicolumn{1}{c|}{\cellcolor[HTML]{ECF4FF}-91396.42} &
  -42035.29
\end{tabular}}
\caption{Deviation-weighted ($F_{dev}$) in the Original and $K$-segment Models.}
\label{tabledyn}
\end{table}

Table~\ref{tabledyn} shows that the $K$-segment models improve fairness by boosting the $F_{dev}$ values when $\alpha>0$, in contrast to the original model.
It can be seen that the absolute values of $F_{dev}$ are relatively high. We apply the $RU_{dev}$ measure, as detailed in Section~\ref{sectiondev} for a clearer discussion.

\end{document}